\documentclass{siamonline1116}

\usepackage{amsfonts}
\usepackage{graphicx}
\usepackage{epstopdf}
\usepackage{mathtools}
\usepackage{subcaption}
\usepackage{algorithm}
\usepackage[noend]{algpseudocode}
\usepackage{amsopn}

\ifpdf
  \DeclareGraphicsExtensions{.eps,.pdf,.png,.jpg}
\else
  \DeclareGraphicsExtensions{.eps}
\fi

\providecommand{\abs}[1]{\left\lvert#1\right\rvert}
\providecommand{\norm}[1]{\left\lVert#1\right\rVert}
\providecommand{\iprod}[1]{\left \langle#1 \right \rangle}
\providecommand{\E}[1]{E\left [#1 \right ]}

\numberwithin{theorem}{section}

\newcommand{\TheTitle}{An Empirical Chaos Expansion Method for Uncertainty Quantification} 
\newcommand{\TheAuthors}{M. Leok and G. Wilkins}

\headers{\TheTitle}{\TheAuthors}

\title{{\TheTitle}\thanks{Submitted to the editors DATE.
\funding{This work was supported in part by NSF grants CMMI-1334759 and DMS-1411792.}}}

\author{
  Melvin Leok\thanks{Department of Mathematics, University of California, San Diego, CA
    (\email{mleok@ucsd.edu}).}
  \and
  Gautam Wilkins\thanks{Department of Mathematics, University of California, San Diego, CA (\email{gwilkins@ucsd.edu}).}
}

\begin{document}
\allowdisplaybreaks

\maketitle

\begin{abstract}
    Uncertainty quantification seeks to provide a quantitative means to understand complex systems that are impacted by parametric uncertainty. The polynomial chaos method is a computational approach to solve stochastic partial differential equations (SPDE) by projecting the solution onto a space of orthogonal polynomials of the stochastic variables and solving for the deterministic coefficients. Polynomial chaos can be more efficient than Monte Carlo methods when the number of stochastic variables is low, and the integration time is not too large. When performing long-term integration, however, achieving accurate solutions often requires the space of polynomial functions to become unacceptably large. This paper presents an alternative approach, where sets of empirical basis functions are constructed by examining the behavior of the solution for fixed values of the random variables. The empirical basis functions are evolved over time, which means that the total number can be kept small, even when performing long-term integration. We introduce this method of empirical chaos expansion, and apply it to a number of model equations, demonstrating that the computational time scales linearly with the final integration time. That is not the case for polynomial chaos in general, since achieving accuracy for long-term integration usually requires larger polynomial bases, causing a nonlinear scaling with the final integration time. We also present an analytical method that uses the dynamics of the SPDE to evolve the empirical basis functions and demonstrate how it can be applied to extend the validity of empirical basis functions in time without the need to sample additional realizations of the original SPDE.
\end{abstract}

\maketitle

\section{Introduction}
Consider the stochastic initial boundary value partial differential equation (SPDE) model problem,
\begin{equation}
u_t(x,t,\xi) = L(u,x,t,\xi)\label{pc_model},
\end{equation}
where $x$ is a spatial variable, $t$ is time, $\xi$ is a random variable with known distribution on probability space $\Omega$, and $L$ is a linear or nonlinear differential operator. SPDEs are used to model systems that contain small scale stochastic components along with large scale deterministic components. Frequently, the deterministic components arise from the governing physics, while the stochastic components are due to measurement errors or some other form of underlying uncertainty. Assuming that the distribution of the stochastic variables is known, we wish to predict the distribution of the solution $u$. SPDEs have frequently demonstrated their use in modeling physical phenomena such as wave propagation \cite{papanicolaou1971wave, kliatskin1980stochastic}, diffusion \cite{pardouxt1980stochastic, papanicolaou1995diffusion, isichenko1992percolation}, Burgers and Navier--Stokes equations with random forcing \cite{bensoussan1973equations, weinan2000statistical, weinan1997probability, weinan2000invariant, mikulevicius2004stochastic, da2003ergodicity, chow1978stochastic, da1994stochastic, sinai1991two, sinai1996burgers}, multivariable predictive control \cite{ni2005multivariable, lou2003estimation, lou2003feedback, lou2004feedback}, and chemical reactors with uncertainties \cite{villegas2012application}.

\subsection{Polynomial Chaos\label{section_pc}}
The original polynomial chaos formulation was introduced by Wiener \cite{wiener1938homogeneous, wiener1966nonlinear}, who used Hermite polynomials to model a Gaussian random process. Polynomial chaos methods begin by choosing a polynomial basis that is orthonormal with respect to the distribution of the random variable $\xi$. Such polynomial bases are known for standard distributions. If we let $\{P_i\}_{i=1}^\infty$ be the orthogonal polynomial basis functions, then we can project the solution $u$ onto the basis functions as follows,
\[
\hat{u}^i(x,t) = \iprod{u,P_i} \coloneqq \int_{\Omega} u(x,t,\xi) P_i(\xi) d\mu,
\]
where $\mu$ is the measure of the probability space $\Omega$. We can then represent the solution $u$ by
\begin{equation}
u(x,t,\xi) = \sum\nolimits_{i=1}^\infty \hat{u}^i(x,t)P_i(\xi)\label{pc_exp}.
\end{equation}
The problem of computing $u(x,t,\xi)$ is equivalent to determining the expansion coefficients $\hat{u}^i(x,t)$. Projecting the initial conditions for $u(x,t,\xi)$ onto the basis functions yield the initial conditions for $\hat{u}^i(x,t)$. Substituting the expansion~\eqref{pc_exp} into the model problem~\eqref{pc_model} yields
\[
\sum\nolimits_{i=1}^\infty \hat{u}_t^i(x,t)P_i(\xi) = L\left(\sum\nolimits_{i=1}^\infty \hat{u}^i(x,t)P_i(\xi), x, t, \xi)\right).
\]
We can then multiply by a test function $P_j$ and take an expectation over $\Omega$ to obtain
\begin{align*}
\hat{u}_t^j(x,t) &= \E{\sum\nolimits_{i=1}^\infty \hat{u}_t^i(x,t)P_i(\xi)P_j(\xi)} = \E{L\left(\sum\nolimits_{i=1}^\infty \hat{u}^i(x,t)P_i(\xi), x, t, \xi)\right)P_j(\xi)},
\end{align*}
where the expectation in the middle simplifies due to the $P_i$'s being orthonormal. We can then truncate the infinite expansion on the right hand side at some finite value $N$, to get
\begin{equation}
\hat{u}_t^j(x,t)  \E{L\left(\sum\nolimits_{i=1}^N \hat{u}^i(x,t)P_i(\xi), x, t, \xi)\right)P_j(\xi)},\quad 1\le j \le N \label{pc_proj}.
\end{equation}
In practice, the right hand side can usually be simplified as well, but it depends on the form of the differential operator $L$. This results in a system of coupled deterministic differential equations for the unknown $\hat{u}^i$'s. Solving the system will give an approximation of the true solution $u$, and this is known as the \textit{stochastic Galerkin method}~\cite{ghanem2003stochastic}.

A result by Cameron and Martin \cite{cameron1944transformations} established that the series \eqref{pc_exp} converges strongly to $L_2$ functionals, which implies that the series converges for stochastic processes with finite second moments. In \cite{ghanem2003stochastic}, Hermite polynomial expansions were coupled with finite element methods to solve SPDEs, and the method was applied to model uncertainty in a variety of physical systems \cite{ghanem1999ingredients, spanos1989stochastic, xiu2003modeling, deb2001solution}, and the convergence properties of polynomial chaos methods for various SPDEs was studied~\cite{lucor2001spectral, chorin1974gaussian, orszag1967dynamical, ernst2012convergence} . An extension by Xiu and Karniadakis \cite{xiu2002wiener} proposed the use of orthogonal polynomials in the Askey scheme \cite{askey1985some} for a variety of other random variable distributions, including uniform and beta distributions. The work also established the same strong convergence results that held for Hermite polynomial expansions of normal random variables. This method is known as generalized polynomial chaos (gPC), and an overview is provided in \cite{Xiu2010}. gPC methods have been applied to a variety of problems in uncertainty quantification, including fluid dynamics and solid mechanics \cite{xiu2002modeling, xiu2002stochastic, zang2002needs, reagan2005quantifying, doostan2007stochastic, lin2006predicting, le2004multi}.

Since their introduction, two well-known issues have become apparent with the gPC method. First, if the number of random variables is too large, then the resulting stochastic Galerkin system becomes too expensive to solve, and Monte Carlo methods outperform gPC methods \cite{wan2005adaptive, gottlieb2008galerkin}. This is the curse of dimensionality, and in this case is due to the need to introduce multidimensional polynomials, one dimension for each random variable. Including even modest orders of high-dimensional polynomials can quickly become prohibitively expensive. Methods have been proposed to limit the growth rate, including sparse truncation \cite{hou2006wiener} that selectively eliminates high-order polynomials, modifications of stochastic collocation methods\cite{doostan2009least}, multi-element gPC expansions \cite{wan2005adaptive, wan2006multi, foo2010multi}, and model reduction to project a high-dimensional random space onto a lower-dimensional one while still preserving important dynamics \cite{doostan2007stochastic, ghanem2007efficient, nouy2007generalized, nouy2008generalized}. While these methods slow the curse of dimensionality, none of them entirely eliminate it, and for systems with a high number of random variables it is more efficient to turn to Monte Carlo methods, which have a convergence rate independent of the dimension of the random space.

The second issue with the gPC method is that in order to accurately approximate a solution with unsteady dynamics over a long time interval, a large number of basis functions must be used, i.e., we must choose a large value for $N$ in \eqref{pc_proj} \cite{gottlieb2008galerkin, najm2009uncertainty, hou2006wiener}. This is not particularly surprising, but it is problematic since the amount of work does not scale linearly with the order of the polynomial basis functions; this issue becomes particularly troublesome as the number of random variables increases. While multi-element gPC expansions \cite{wan2005adaptive, wan2006multi} have been suggested as an option, they do not entirely solve this issue. Another approach is the time-dependent polynomial chaos method \cite{gerritsma2010time, heuveline2011local}, which integrates for a short period of time, then treats the solution as a new random variable and attempts to numerically derive an orthogonal polynomial basis for the new random variable. Later, we will introduce a method that addresses this issue by combining model reduction techniques with gPC, and discuss a basis evolution technique that further reduces the computational cost of the approach.

\subsection{General Basis}
The above results can be extended in a straightforward manner to non-orthonormal bases as well. Consider any set of functions $\{\Psi^i\}_{i=1}^\infty$ that form a basis for $L^2(\Omega)$, but are not necessarily orthonormal. Then, there exist coefficients $\hat{u}^i(x,t)$ such that
\[
u(x,t,\xi) = \sum\nolimits_{i=1}^\infty \hat{u}^i(x,t) \Psi^i(\xi),
\]
which implies that
\begin{align*}
\iprod{u, \Psi^i} &= \iprod{\sum\nolimits_{j=1}^\infty \hat{u}^j\Psi^j, \Psi^i}= \sum\nolimits_{j=1}^\infty \hat{u}^j\iprod{\Psi^j, \Psi^i}.
\end{align*}
If we let
\[
A_{ji} = \iprod{\Psi^j, \Psi^i},\quad \hat{u} = \begin{pmatrix}\hat{u}^1\\\hat{u}^2\\\vdots\\\end{pmatrix},\quad \mathrm{and}\quad f = \begin{pmatrix}\iprod{u, \Psi^1}\\\iprod{u, \Psi^2}\\\vdots\end{pmatrix},
\]
then we can write the above system as
\begin{equation}
A \hat{u} = f\label{gen_ic},
\end{equation}
where the $\hat{u}$ vector is unknown, and the $A$ and $f$ matrices can be computed by inner products (if $u$ is known). Since the model problem is an initial boundary value problem, we can compute the initial values for the $\hat{u}_i$'s. Just like with standard polynomial chaos, we can derive a stochastic Galerkin method by substituting the expansion into the model problem~\eqref{pc_model}, multiplying by a test function $\Psi^j$, and taking an expectation over $\Omega$. If we truncate the expansion at finite $N$, we obtain
\begin{equation}
\E{\sum\nolimits_{i=1}^N \hat{u}_t^i(x,t) \Psi^i(\xi)\Psi^j(\xi)} = \E{L\left( \sum\nolimits_{i=1}^N \hat{u}^i(x,t) \Psi^i(\xi),x,t,\xi \right)\Psi^j(\xi)}\label{gen_exp}.
\end{equation}
Now, if we let
\[
\tilde{A}_{ji} = \iprod{\Psi^j, \Psi^i},\quad \hat{u} = \begin{pmatrix}\hat{u}^1\\\hat{u}^2\\\vdots\\\hat{u}^n\\\end{pmatrix},\quad \mathrm{and}\quad b_j = \E{L\left( \sum\nolimits_{i=1}^N \hat{u}^i(x,t) \Psi^i(\xi),x,t,\xi \right)\Psi^j(\xi)},
\]
then, we obtain the system
\begin{equation}
\tilde{A}\hat{u}_t = b,
\label{sg_sys}
\end{equation}
which is a deterministic implicit differential equation whose solution is the $\hat{u}^i$ coefficients.

\subsection{Conversion Between Bases}
\label{basis_conv}
Assume we have two distinct bases $\{\Psi^i\}_{i=1}^\infty$ and $\{\Phi^i\}_{i=1}^\infty$ for $L^2(\Omega)$, such that,
\begin{align*}
u(x,t,\xi) &= \sum\nolimits_{i=1}^\infty \hat{u}^i_{\Phi}(x,t)\Phi^i(\xi),\qquad
u(x,t,\xi) = \sum\nolimits_{i=1}^\infty \hat{u}^i_{\Psi}(x,t)\Psi^i(\xi).
\end{align*}
If we know $\hat{u}^i_{\Psi}$, we can compute $\hat{u}^i_{\Phi}$ with a standard change of basis operation,
\begin{align*}
\iprod{u, \Phi^j} &= \sum\nolimits_{i=1}^\infty \iprod{\hat{u}^i_\Psi(x,t)\Psi^i, \Phi^j}= \sum\nolimits_{i=1}^\infty \hat{u}^i_\Psi(x,t)\iprod{\Psi^i, \Phi^j}.
\end{align*}
Letting
\[
M_{ji} = \iprod{\Psi^i, \Phi^j},\quad\mathrm{and}\quad b_j = \iprod{u, \Phi^j},
\]
implies that
\begin{equation}
b = M\hat{u}_{\Psi}\label{change_basis}.
\end{equation}
The right hand side is a straightforward computation, since $\hat{u}_{\Psi}$ is known, and we can use \eqref{gen_ic} to compute $\hat{u}_{\Phi}$. The same computation applies when the two bases are finite-dimensional instead of infinite-dimensional. In the finite-dimensional case, the operations above first project $u$ onto the subspace spanned by $\{\Psi^j\}$, followed by a projection onto the subspace spanned by $\{\Phi^i\}$.

\section{Empirical Bases\label{emp_section}}
We now seek to resolve the issue of long-term integration for gPC methods, using the following general approach. First, we come up with a finite and small set of basis functions that result in a low projection error over a short time interval (say $0\le t \le \tau_0$). Then, perform a stochastic Galerkin method via \eqref{gen_exp} in order to approximate the solution up to time $\tau_0$. This will give values for the $\hat{u}^i$ coefficients. Now, come up with a different set of basis functions that result in a low projection error over the short time interval $\tau_0 \le t \le \tau_1$. We can use the previous solution at $\tau_0$ as the initial condition, and calculate the initial values of the $\hat{u}^i$'s in the new basis by applying the change of basis equation~\eqref{change_basis}. We can then perform another stochastic Galerkin method with the new basis to compute the solution up to time $\tau_1$. Continuing this process iteratively allows us to compute the solution out to a long time $t$ with a low error, and at each individual step, we only require a relatively small basis.

The principal issue is determining a subspace spanned by a small number of basis functions for each interval $\tau_i \le t \le \tau_{i+1}$, such that the solution $u$ will have a low projection error in the subspace. Once we have such a set of basis functions for each time interval, then we can follow the approach outlined above to solve the SPDE on a relatively small set of basis functions.

A similar attempt is made in \cite{gerritsma2010time}, but it relies on integrating the solution to a small time in the future, and then treating the solution as a new random variable that must be added to the gPC expansion before integrating out to a later time. This entails numerically deriving a new set of orthogonal polynomial basis functions that are inferred from the numerical probability distribution of the solution. In this paper, we present an alternate approach, where we no longer require the basis functions to be orthonormal, and instead use techniques of model reduction to construct an empirically optimal set of basis functions at each given timestep.

\subsection{Sampling Trajectories}
\label{sample_trajectory_sec}
Fix a single point $(x_0,t_0)$ in spacetime, and examine how the function $u(x_0, t_0, \xi)$ varies in the random variable $\xi$. If we examine many such sample trajectories, then we can choose a good subspace for projecting the solution by requiring that it have a low projection error for most or all of the sample trajectories.

We can find such trajectories by choosing a set of fixed values for $\xi$, say $\{\xi_l\}_{l=1}^K$, and then solving the resulting deterministic PDE for each $\xi_l$, much like in Monte Carlo methods. The difference is that Monte Carlo methods require solving the PDE for an extremely large number of values of $\xi$ before they converge, and we use much fewer sample trajectories. Once we know the solutions for each of the values of $\{\xi_l\}$, we can construct the following matrix,
\begin{equation}
T=
\begin{pmatrix}
u(x_1,t_1,\xi_1) & u(x_1,t_1,\xi_2) & \ldots & u(x_1,t_1,\xi_K)\\
u(x_1,t_2,\xi_1) & u(x_1,t_2,\xi_2) & \ldots & u(x_1,t_2,\xi_K)\\
\vdots & \vdots & & \vdots\\
u(x_1,t_N,\xi_1) & u(x_1,t_N,\xi_2) & \ldots & u(x_1,t_N,\xi_K)\\
u(x_2,t_1,\xi_1) & u(x_2,t_1,\xi_2) & \ldots & u(x_2,t_1,\xi_K)\\
\vdots & \vdots & & \vdots\\
u(x_M,t_N,\xi_1) & u(x_M,t_N,\xi_2) & \ldots & u(x_M,t_N,\xi_K)
\end{pmatrix},\label{pod_sys}
\end{equation}
where we have discretized spacetime into the set of points $\{(x_i, t_j)\}_{(i,j)=(1,1)}^{(M,N)}$. Each row of the matrix $T$ represents a discrete trajectory of the solution as a function of the random variable $\xi$. Note that, in general, this matrix will have far more rows than columns. We now seek a, hopefully small, set of basis functions in $\xi$ that will result in a low projection error for each of the rows of $T$. The motivating assumption is that such a set of basis functions will also result in a low projection error for the true solution $u(x,t,\xi)$.

\subsection{Proper Orthogonal Decomposition}
\label{pod}
Proper Orthogonal Decomposition (POD) is a standard tool in the field of model reduction \cite{chatterjee2000introduction, kerschen2005method} and is also known as the Karhunen--Lo\`eve transform \cite{karhunen1947lineare, kosambi1943statistics}, or Principal Component Analysis \cite{person1901lines} in the finite-dimensional case. It is often used to construct low-dimensional representations of very high-dimensional spaces that preserve the most important dynamics \cite{liang2002proper, liang2002proper2, willcox2002balanced}.

Following the exposition from \cite{liang2002proper}, in the discrete case we have $T$ in $\mathbb{R}^{NM} \times \mathbb{R}^K$ from \eqref{pod_sys}, where we view each row of $T$ as a sample from a random vector that takes values in $\mathbb{R}^K$. We then seek a set of $k<K$ basis vectors such that the rows of $T$ can be projected with minimum mean square error onto the $k$ basis vectors. The optimization problem can be expressed as
\[
\min_{P_k}\ \E{\norm{v - P_kv}^2},
\]
where $v$ is a random vector from the space sampled by the rows of $T$, the expectation is over the rows of $T$, $P_kv$ is the projection of $v$ onto the $k$ basis vectors, and the minimization is over all sets of $k$ basis vectors.  In practice, we perform a singular value decomposition (SVD),
\[
T = U\Sigma V^T.
\]
Then, the first $k$ columns of the right singular matrix $V$ are the optimal set of $k$ basis vectors that minimize the mean square error when we project the $NM$ rows of $T$ onto their span. This will give us a set of $k$ vectors that span a subspace that will result in a low projection error for all of the rows of the $T$ matrix. We can now construct a basis function in $\xi$ for each row of the $V$ matrix by constructing an interpolant through its values at each of the $\xi_l$ points.

\subsection{Empirical Chaos Expansion}
At this point, we can construct a general algorithm to solve the SPDE model problem~\eqref{pc_model} using successive sets of empirical basis functions. We refer to this method as \textit{empirical chaos expansion}. We begin by choosing a set, $\{\xi_l\}_{l=1}^K$, from the range of the random variable $\xi$. If there were multiple random variables, we would choose a set of random vectors from the product space of the random variables. If we wish to solve the SPDE from time $t_0$ to time $t_f$, we partition the interval  by choosing nodal points $t_0<\tau_1<\tau_2<\cdots<\tau_{n-1}<t_f$. We then solve the deterministic partial differential equation for each fixed value of $\xi$ in the set $\{\xi_l\}_{l=1}^K$, starting at time $t_0$ and ending at time $\tau_1$. Assuming that the solutions are computed using a finite-difference method, we will have values for the solution $u(x,t,\xi)$ at a set of points $(x_i, t_j, \xi_l)$. Fixing a point in spacetime and varying the value of $\xi$ yields a row of the matrix $T$ in Section~\ref{sample_trajectory_sec}, and we use the full set of solutions to construct the full $T$ matrix. We then perform a POD of the $T$ matrix by computing its singular value decomposition, and choosing the first $N_1$ columns of the right singular value matrix to be the set of basis functions. In practice, $N_1$ is chosen so that the $(N_1+1)$-st singular value is small relative to the largest singular value. Each column of the right singular value matrix can be viewed as a table of values for the basis function at each point in  $\{\xi_l\}_{l=1}^K$.

We could construct an interpolant for each of the basis column vectors, but it is computationally faster to use a numerical quadrature formula with $\{\xi_l\}_{l=1}^K$ as the quadrature points. We denote the set of empirical basis functions for the time interval $[t_0, \tau_1]$ by $\{\Psi^i_1(\xi)\}_{i=1}^{N_1}$. For the first time interval $[t_0, \tau_1]$, we can approximate the true solution $u$ by $u(x,t,\xi) \approx \sum\nolimits_{i=1}^{N_1} \hat{u}^i_1(x,t) \Psi^i_1(\xi)$. To determine the initial values of $\hat{u}^i_1$, we project the initial condition onto each of the basis functions, $\hat{u}^i_1(x,t_0) = \iprod{u(x,t_0), \Psi^i_1}$, and we solve the propagation equation~\eqref{sg_sys} for the stochastic Galerkin system for the SPDE, which computes $\hat{u}^i_1$ up to time $\tau_1$. The expectations in the propagation equation can be computed using numerical quadrature.

The solution at the final time $\tau_1$ can be used as the initial condition on the interval $[\tau_1, \tau_2]$. We begin by computing the solution to the deterministic partial differential equation for each fixed value of $\xi$ in the set $\{\xi_l\}_{l=1}^K$, starting at time $\tau_1$ and ending at time $\tau_2$. We now use the computed solutions to construct a new $T$ matrix that has trajectories over the time interval $[\tau_1, \tau_2]$. We then perform a POD of the $T$ matrix by computing its singular value decomposition, and choose the first $N_2$ columns of the right singular value matrix to be the set of basis functions. This gives a set of empirical basis functions $\{\Psi^i_2(\xi)\}_{i=1}^{N_2}$. Although we have values $\hat{u}^i_1(x,t)$ up to time $\tau_1$, these are the coefficients in the old basis $\{\Psi^i_1(\xi)\}_{i=1}^{N_1}$. We can convert them to coefficients in the new basis $\{\Psi^i_2(\xi)\}_{i=1}^{N_2}$ by following the method described in Section~\ref{basis_conv}. Recall that applying the method requires us to compute entries of $M_{ji} = \iprod{\Psi^i_1, \Psi^j_2}$. Again, these inner products can be computed using numerical quadrature. Then, the initial $\hat{u}_2$ coefficient vector in the new basis is related to the coefficient vector $\hat{u}_1$ in the old basis by $\hat{u}_2=M \hat{u}_1$. Now, the true solution $u$ is approximated by $u(x,t,\xi) \approx \sum\nolimits_{i=1}^{N_2} \hat{u}^i_2(x,t) \Psi^i_2(\xi)$, on the interval $[\tau_1, \tau_2]$. Since we only have $\hat{u}^i_2$ at time $\tau_1$, we solve the propagation equation~\eqref{sg_sys}, which computes the $\hat{u}^i_2$ coefficients up to time $\tau_2$.

At this point we repeat the procedure in the previous paragraph to generate a new set of empirical basis functions for the interval $[\tau_2, \tau_3]$, project the old coefficients onto the new basis, and solve the new propagation equation to obtain a solution up to time $\tau_3$. This process is then repeated for every timestep until we reach the final integration time.

\subsection{Convergence}
Let $u_m^k = \sum\nolimits_{i=1}^m \hat{u}^i(x,t) \Psi^i(\xi)$ be the approximation to the true solution on a time interval $[\tau_{i}, \tau_{i+1}]$ constructed through empirical chaos by sampling a discrete set of $k$ trajectories in the space of random variables, $\{\xi_l\}_{l=1}^K$. Let $\bar{u}_m^k$ be the projection of the true solution $u$ onto the span of the $k$ trajectories. Note that the projection error at any of the $\xi_l$ points will be zero. Then, from the triangle inequality we have
\[
\norm{u-u^k_m} \le \norm{u-\bar{u}^k_m} + \norm{\bar{u}^k_m - u^k_m},
\]
in some appropriate norm. The error in the first term can be controlled by choosing $\{\xi_l\}_{l=1}^K$ as interpolation nodes and increasing $K$. This follows from standard interpolation theory, assuming sufficient regularity of the solution $u$. The error in the second term is due to POD, where instead of projecting $u$ onto the full span of the trajectories, we project onto a $m$-dimensional subspace. This error can be controlled by increasing $m$.

\subsection{Computational Cost}
In order to fully solve the SPDE at an intermediate timestep $j+1$ (i.e. any timestep other than the first), we must perform the following operations:
\begin{enumerate}
    \item Solve the deterministic PDE for each point in $\{\xi_l\}_{l=1}^K$. This operation is simply $K$ times the cost of solving a single instance of the deterministic PDE. However, this is trivially parallelizable, which can dramatically reduce the computational time.
    \item Compute the POD of the $T$ matrix by computing its SVD. The $T$ matrix always has $K$ columns, but the number of rows depends on the number of finite-difference grid points and is almost always much larger than $K$. Since we only need the right singular values, which is a $K\times K$ matrix, we can reduce the cost of the SVD operation by only computing the singular value matrix $\Sigma$, and the right singular value matrix $V$.
    \item Project the old $\hat{u}^i_j$ coefficients onto the new set of basis functions. This requires us to compute the $M$ matrix, which will be $N_{j+1}\times N_{j}$ (where $N_{j}$ and $N_{j+1}$ are the number of empirical basis functions for the $j$-th and $(j+1)$-st timesteps, respectively). Each entry of the matrix is computed by calculating an inner product with a numerical quadrature method. We then multiply the $M$ matrix by the old $\hat{u}^i_j$ coefficients.
    \item Evaluate the propagation equation~\eqref{sg_sys}. This will require computing expectations, which can be done with numerical quadrature. The cost of this operation depends on the number of basis functions and also how complicated the differential operator $L$ is.
    \item Solve the propagation equation~\eqref{sg_sys} to advance the solution to the next timestep. The cost of this operation is highly dependent on the differential operator $L$. In some cases, this is no harder than solving the deterministic PDE, and in others it is far more computationally expensive. In either case, it is a coupled system of $N_{j+1}$ PDEs.
\end{enumerate}
The key advantage of this approach is that we can control the number of empirical basis functions by adjusting the size of the timestep. gPC uses the same set of basis functions for every point in time, and for many practical SPDEs this means that achieving accuracy in the solution out to long integration times requires using a very large set of basis functions. In contrast, each set of empirical basis functions only needs to have a low projection error over its corresponding time interval, which can be quite short. This allows the empirical chaos expansion method to run in linear time as a function of the final integration time, i.e., doubling the final integration time doubles the computational cost of empirical chaos expansion.

\subsection{Computing Solution Statistics}
If there are $n$ total timesteps then we will need to keep track of a set of $n$ sets of basis functions and their associated coefficients in order to fully reconstruct the solution. In the end, we still have a fully analytic representation of the true solution $u$, and we can use it to compute any desired solution statistics. The actual computations will be more complex than those for gPC, due to the fact that the basis functions are not orthogonal, and there are multiple sets of basis functions (one for each timestep). Let $\{\Psi^i_j\}_{i=1}^{N_j}$ be the set of $N_j$ basis functions for timestep $j$, and let $\{\hat{u}^i_j\}_{i=1}^{N_j}$ be the corresponding set of coefficients. If we want to compute the mean of the true solution, then we compute,
\begin{align*}
    \E{u(x,t,\xi)} &= \E{\sum\nolimits_{i=1}^{N_{k^*}} \hat{u}^i_{k^*} \Psi^i_{k^*}}= \sum\nolimits_{i=1}^{N_{k^*}} \hat{u}^i_{k^*} \E{\Psi^i_{k^*}},
\end{align*}
where $k^*$ is the timestep that contains the time $t$, and $\E{\Psi^i_{k^*}}$ can be computed using numerical quadrature. Other solution statistics can be computed in a similar manner.

\section{One-Dimensional Wave Equation\label{owsection}}
Consider the SPDE,
\begin{align}
u_t(x,t,\xi) &= \xi u_x(x,t,\xi),\quad 0\le x\le 2\pi,\quad t\ge 0,\label{ow_wave}\\
u(x,0,\xi) &= \cos(x),\label{ow_wave_ic}
\end{align}
with periodic boundary conditions, and $\xi$ is uniform on $[-1,1]$. The exact solution is
\begin{equation}
u(x,t,\xi)=\cos(x-\xi t).\label{ow_wave_exact}
\end{equation}
Note that over the space for this problem,
\begin{align*}
\iprod{f,g} &\coloneqq \int_{-1}^1 f(\xi)g(\xi)(1/2)d\xi, \qquad \E{f} \coloneqq \int_{-1}^1 f(\xi)(1/2)d\xi.
\end{align*}
This hyperbolic SPDE was analyzed in \cite{gottlieb2008galerkin}, and one of the findings was that for a fixed polynomial basis, the error scales linearly with final integration time, thus requiring the gPC expansion to continually add additional terms in order to achieve accurate long-term solutions.

\subsection{Polynomial Chaos Expansion}
\label{wave_gpc}
\begin{figure}[htb]
	\centering
	\begin{subfigure}[b]{0.325\textwidth}
		\centering
		\includegraphics[scale=0.25]{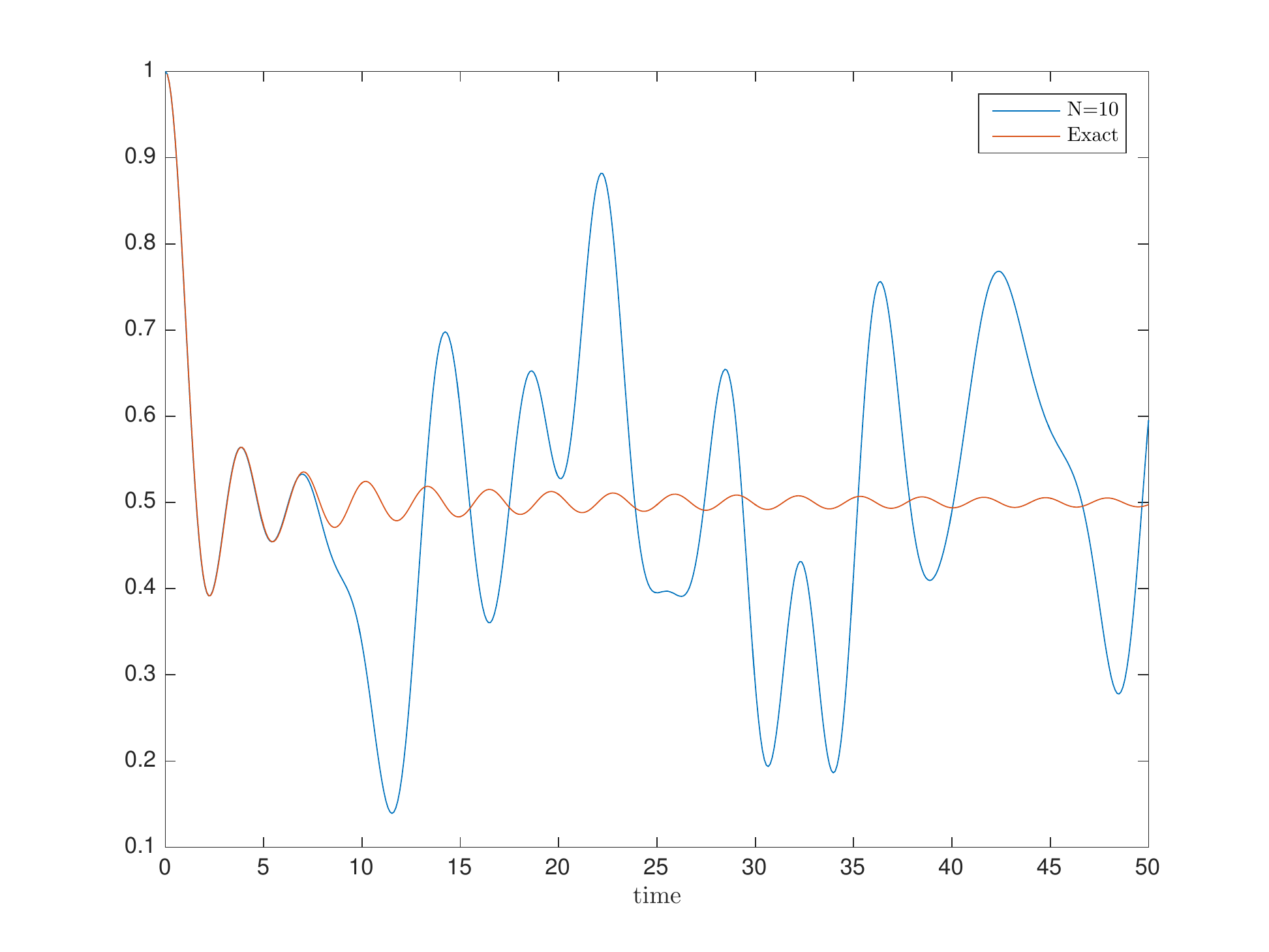}
		\caption{\label{mse_pc_plot_10}$d=10$}
	\end{subfigure}
	\begin{subfigure}[b]{0.325\textwidth}
  		\centering
  		\includegraphics[scale=0.25]{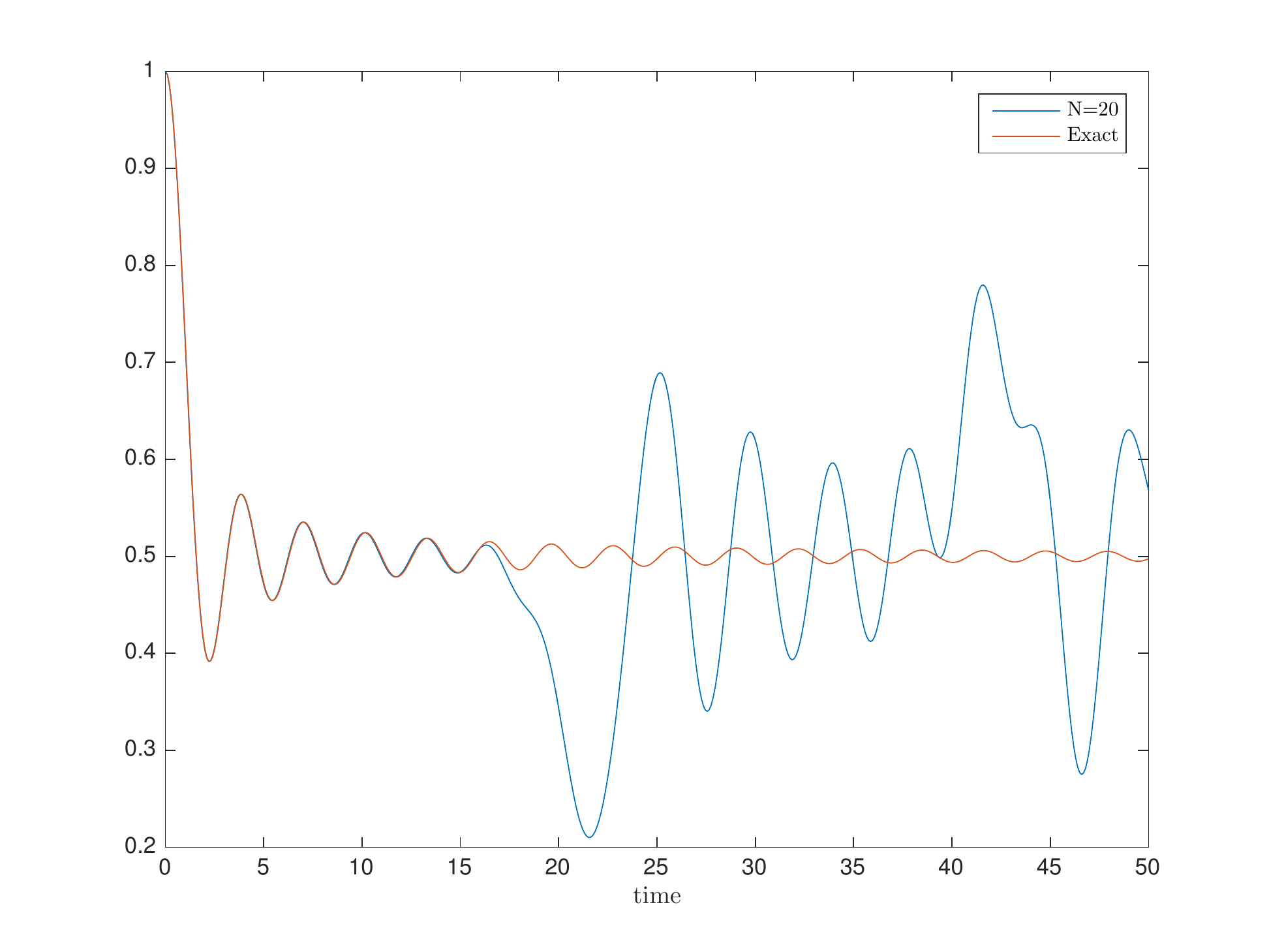}
  		\caption{\label{mse_pc_plot_20}$d=20$}
	\end{subfigure}
	\begin{subfigure}[b]{0.325\textwidth}
		\centering
  		\includegraphics[scale=0.25]{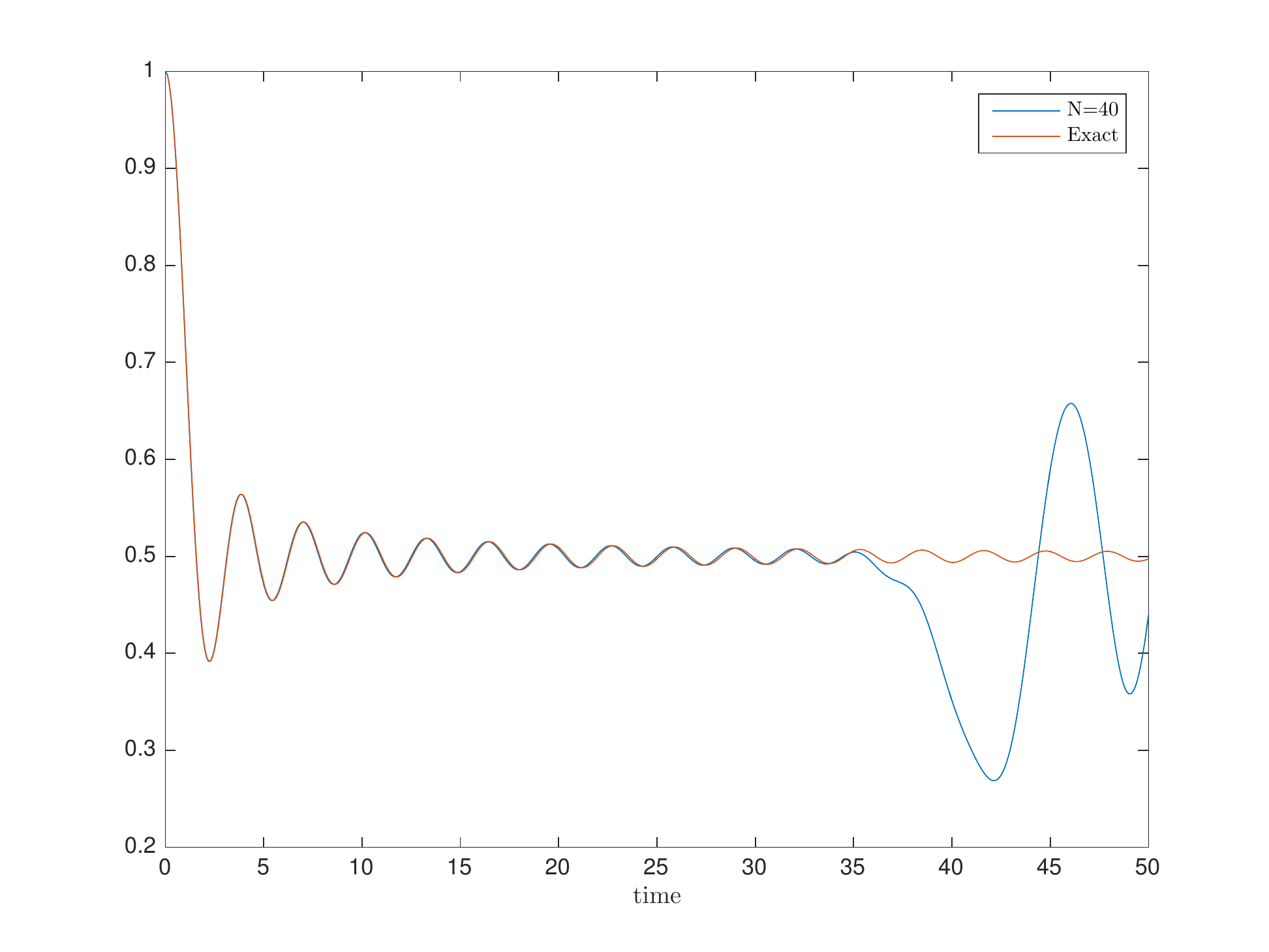}
  		\caption{\label{mse_pc_plot_40}$d=40$}
  	\end{subfigure}
    \caption{Mean Square Expectation at $x=0$ of solution to \eqref{ow_wave}, computed using gPC with stochastic Galerkin using Legendre polynomials up to order $d$.}  
\end{figure}

When using gPC, we choose the normalized Legendre polynomials as our basis functions. Let $L^i$ be the $i$-th normalized Legendre polynomial and apply stochastic Galerkin as in Section~\ref{section_pc}, and expand the true solution as
\[
u(x,t,\xi) = \sum\nolimits_{i=1}^\infty \hat{u}^i(x,t)L^i(\xi).
\]
Substituting into \eqref{ow_wave}, multiplying by a test function $L^j(\xi)$, and taking the expectation yields
\begin{align*}
 \hat{u}^j_t &= \sum\nolimits_{i=1}^\infty \hat{u}^i_t(x,t) \E{L^i(\xi)L^j(\xi)} = \E{\sum\nolimits_{i=1}^\infty \hat{u}^i_t(x,t)L^i(\xi)L^j(\xi)}\\
 &= \E{\xi \sum\nolimits_{i=1}^\infty \hat{u}^i_x(x,t)L^i(\xi)L^j(\xi)} =
 \sum\nolimits_{i=1}^\infty \hat{u}^i_x(x,t) \E{\xi L^i(\xi)L^j(\xi)}.
\end{align*}
Letting
\[
A_{ji} = \E{\xi L_jL_i}, \quad \mathrm{and} \quad \hat{u} =
\begin{pmatrix}
    \hat{u}^1(x,t)\\
    \hat{u}^2(x,t)\\
    \vdots\\
\end{pmatrix},
\]
simplifies the previous system to
\begin{equation}
\hat{u}_t = A\hat{u}_x\label{ow_wave_pc}.
\end{equation}
Since the initial condition is deterministic, $\hat{u}^1(x,0) = \cos(x)$, and $\hat{u}^i(x,0)=0$ for $i>1$. We truncate the infinite system to finite $N$ and solve the resulting deterministic system of PDEs. The exact mean square expectation at $x=0$ is given by
\begin{align}
\E{(u(0,t,\xi)^2} &= \int_{-1}^1 (\cos^2(\xi t))(1/2)d\xi= \frac{1}{2}\left(1+\frac{\cos(t)\sin(t)}{t} \right)\label{mse_exact}.
\end{align}
See Figures~\ref{mse_pc_plot_10}, \ref{mse_pc_plot_20}, \ref{mse_pc_plot_40} for a comparison of the exact value of the mean square expectation with the numerical value computed by the stochastic Galerkin method with different values of $N$. Past a certain point in time all of the PC expansions diverge from the exact value.

\begin{figure}[htb]
	\centering
    \includegraphics[scale=0.35]{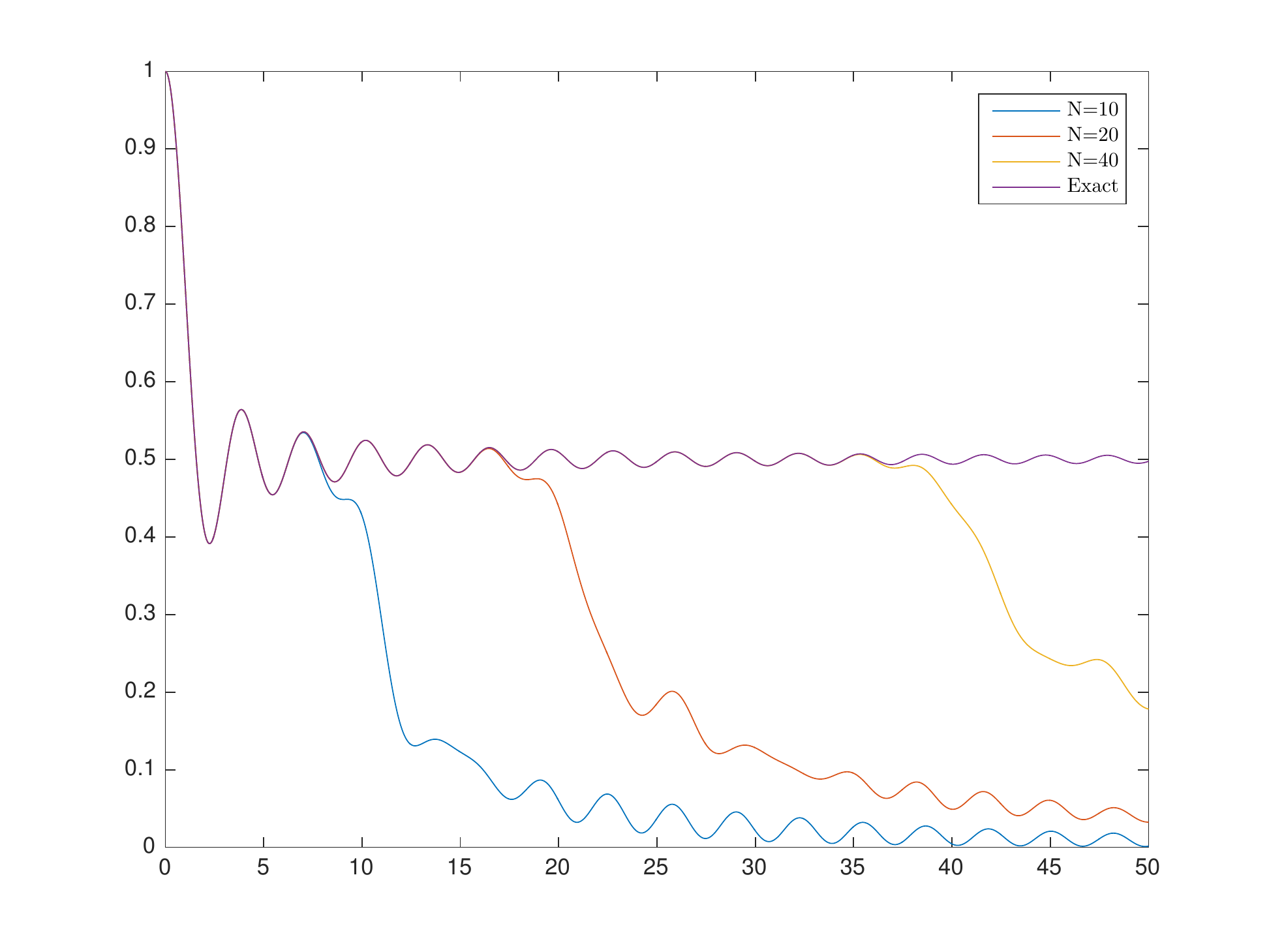}
    \caption{\label{mse_exact_plot}Mean Square Expectation at $x=0$ for \eqref{ow_wave}. Projection of the exact solution onto gPC basis functions for increasing polynomial order $N$.}
\end{figure}

This is not an issue with the stochastic Galerkin method. In fact, even if we project the exact solution \eqref{ow_wave_exact} onto the space of Legendre polynomials (see Figure~\ref{mse_exact_plot}) and compute its mean square expectation, we can see that more and more polynomial basis functions are needed to accurately project the exact solution. 

This raises another important concern with standard polynomial chaos techniques. If the initial condition is stochastic, then we may need a large number of polynomial functions \emph{just to accurately project the initial condition}. For example, if we chose the initial condition to be the exact solution $\cos(x-\xi t)$ for a large value of $t$, then we would need a large number of Legendre polynomials even to accurately compute the solution over a short time interval.

\subsection{Empirical Chaos Expansion}
\label{wave_emp}
In order to apply an empirical chaos expansion, we must account for the fact that the basis functions might not be orthogonal. This results in the stochastic Galerkin system becoming slightly more complicated than the one for gPC. If we let $\{\Psi^i(\xi)\}_{i=1}^N$ be a set of possibly non-orthogonal basis functions, the exact solution to \eqref{ow_wave} can be approximated by 
\begin{equation}\label{empirical_chaos_expansion}
	u(x,t,\xi) \approx \sum\nolimits_{i=1}^N \hat{u}^i(x,t)\Psi^i(\xi).	
\end{equation}
Substituting this into \eqref{ow_wave} results in $\sum\nolimits_{i=1}^N \hat{u}^i_t(x,t)\Psi^i(\xi) = \xi \sum\nolimits_{i=1}^N \hat{u}^i_x(x,t)\Psi^i(\xi)$.
Multiplying by a test function $\Psi^j(\xi)$, and taking an expectation over $\Omega$ yields
\begin{multline}
\sum\nolimits_{i=1}^N \hat{u}^i_t(x,t)\E{\Psi^i(\xi)\Psi^j(\xi)} = \E{\sum\nolimits_{i=1}^N \hat{u}^i_t(x,t)\Psi^i(\xi)\Psi^j(\xi)}\\= \E{\sum\nolimits_{i=1}^N \hat{u}^i_x(x,t)\Psi^i(\xi)\Psi^j(\xi)\xi}= \sum\nolimits_{i=1}^N \hat{u}^i_x(x,t)\E{\Psi^i(\xi)\Psi^j(\xi)\xi}\label{ow_wave_emp_exp}.
\end{multline}
Letting
\[
M_{ji} = \E{\Psi^i(\xi)\Psi^j(\xi)},\quad A_{ji} = \E{\Psi^i(\xi)\Psi^j(\xi)\xi}, \quad \mathrm{and} \quad
u =
\begin{pmatrix}
u^1(x,t)\\
u^2(x,t)\\
\vdots\\
u^N(x,t)
\end{pmatrix},
\]
then \eqref{ow_wave_emp_exp} can be expressed as
\begin{equation}
    Mu_t = Au_x.\label{ow_wave_sys}
\end{equation}
If we know the empirical basis functions $\Psi^i$ and the initial value of $\hat{u}^i$, then the matrices $M$ and $A$ can both be computed, and \eqref{ow_wave_sys} can be solved to obtain the coefficients at later times.

We follow the method in Section~\ref{emp_section} to construct empirical basis functions over small time intervals by sampling trajectories and applying POD. Fixing $\xi$ makes the PDE deterministic, which can be solved with existing numerical methods. We use a method of lines discretization to approximate the spatial derivative and apply a fourth-order Runge--Kutta method to perform the time integration (\texttt{ode45} in MATLAB). The POD is truncated, and the empirical basis functions $\{\Psi^i(\xi)\}$ are generated. We do not interpolate the basis functions, but instead use a composite trapezoidal quadrature rule to compute the expectations and inner products.

\begin{figure}[htb]
	\centering
	\begin{subfigure}[b]{0.48\textwidth}
  		\centering
    	\includegraphics[scale=0.35]{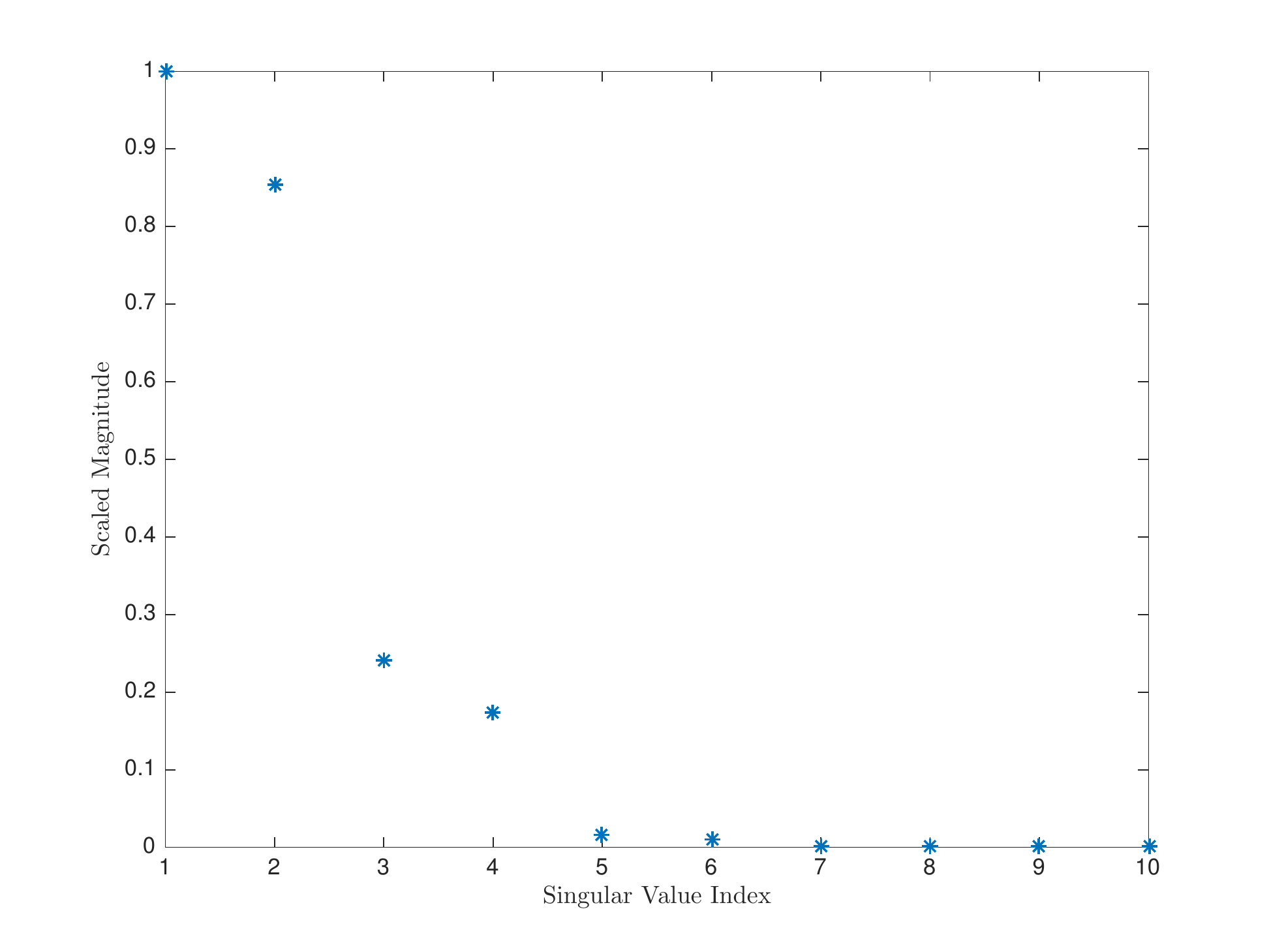}
    	\caption{\label{ec_svs_ts1}Scaled singular values from POD for the first timestep.}
	\end{subfigure}\quad
	\begin{subfigure}[b]{0.48\textwidth}
  		\centering
    	\includegraphics[scale=0.35]{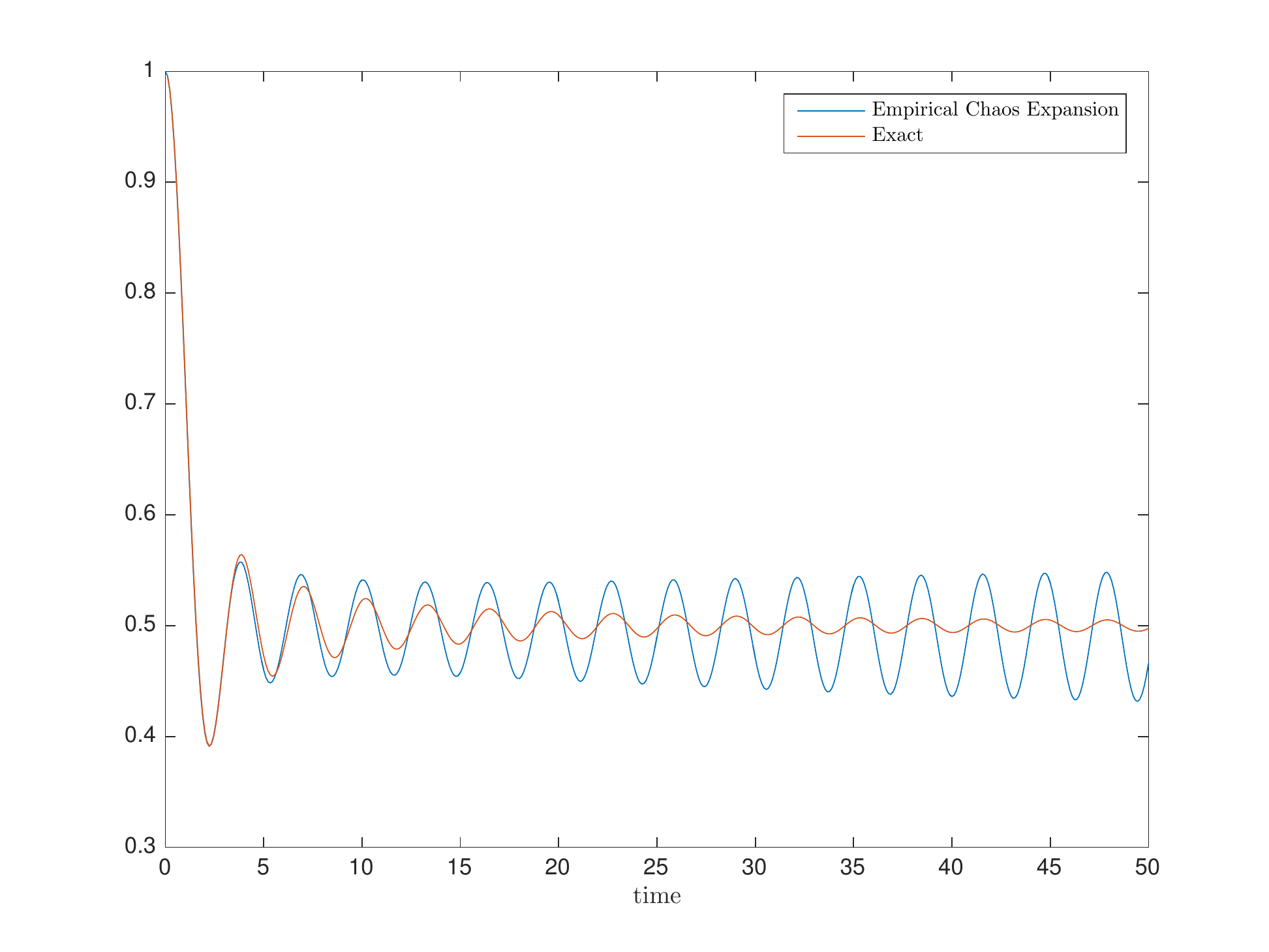}
    	\caption{\label{mse_ec_ts1_tf50_nb_4}Maximum of 4 empirical basis functions on each timestep.}
	\end{subfigure}\\
	\begin{subfigure}[b]{0.48\textwidth}
		\centering
    	\includegraphics[scale=0.35]{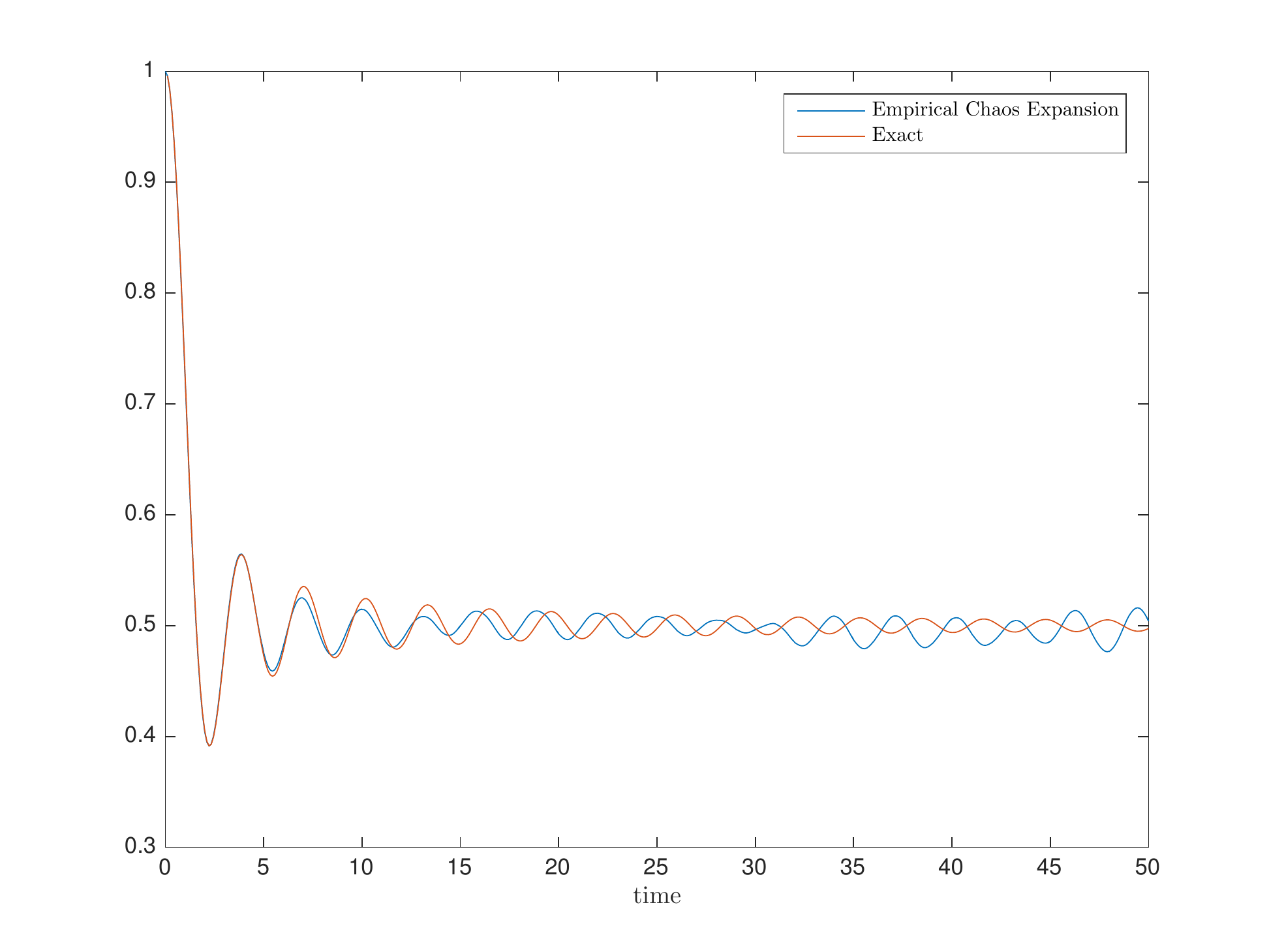}
    	\caption{\label{mse_ec_ts1_tf50_nb_5}Maximum of 5 empirical basis functions on each timestep.}
    \end{subfigure}\quad
 	\begin{subfigure}[b]{0.48\textwidth}   
 		\centering
		\includegraphics[scale=0.35]{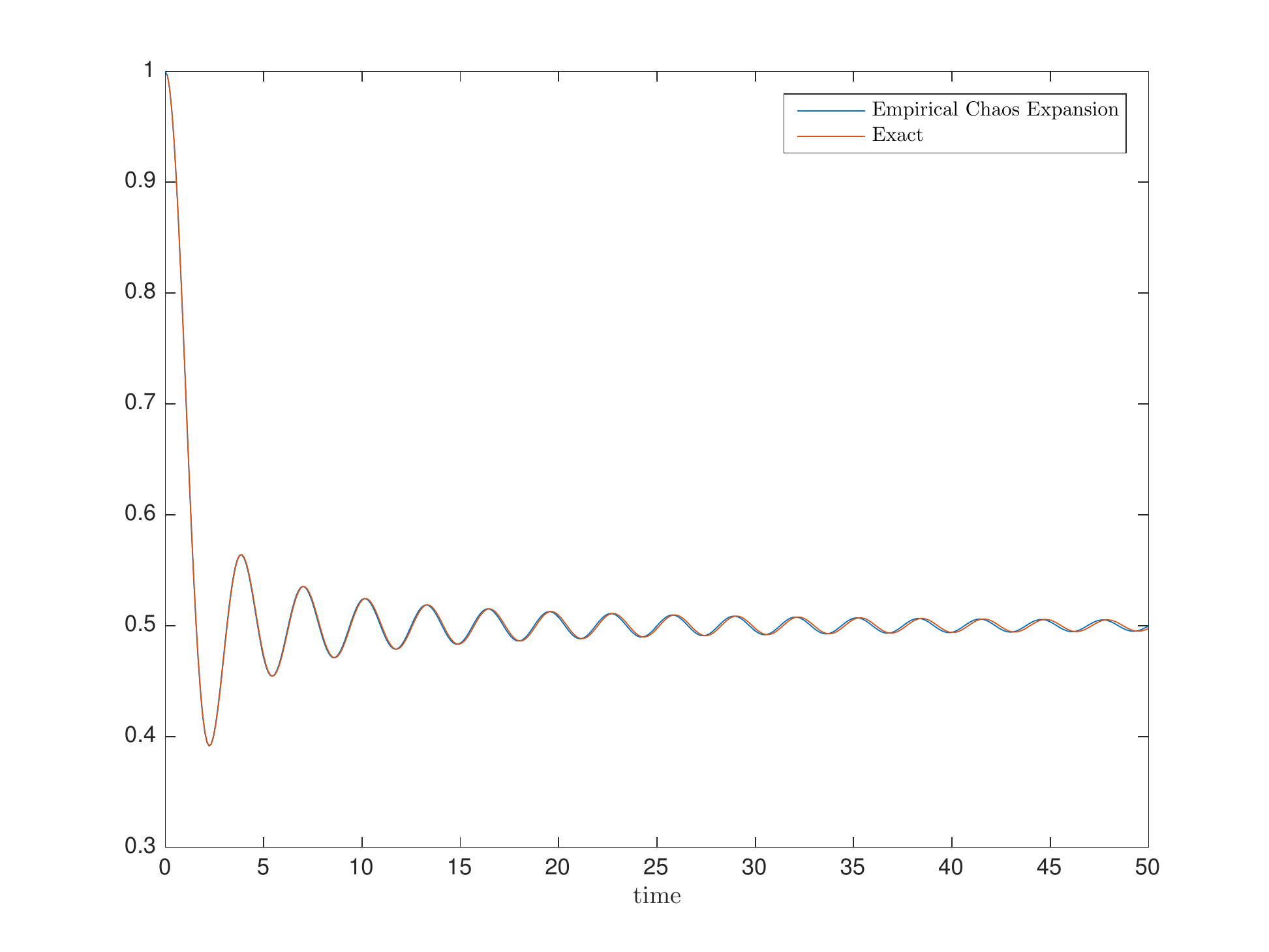}
    	\caption{\label{gc_ts1_tf50_plot}Maximum of 9 empirical basis functions on each timestep.}
	\end{subfigure}
    \caption{Mean Square Expectation at $x=0$ of the solution to \eqref{ow_wave}, computed using empirical chaos expansion with stochastic Galerkin.}
\end{figure}

To determine an appropriate number of basis functions to use for a given timestep, we examine the singular values from POD in Figure~\ref{ec_svs_ts1}, scaled so that the first singular value is 1. There is a sharp drop-off in the singular values around the fifth singular value. In this case, we choose to truncate the expansion when the scaled singular values drops below $10^{-4}$, which corresponds to the ninth singular value. The threshold is arbitrary and was determined empirically by truncating the expansion at varying values. This could possibly be made more precise by recalling that the squared Frobenius norm of the difference between a matrix and its truncated SVD is given by the sum of squares of the ignored singular values.

Figures~\ref{mse_ec_ts1_tf50_nb_4}, \ref{mse_ec_ts1_tf50_nb_5}, and \ref{gc_ts1_tf50_plot} show the result of applying empirical chaos expansion to the one-dimensional wave equation with 4, 5, and 9 basis functions, respectively. These used 120 Chebyshev nodes for $\xi$ in the interval $[-1,1]$ to compute trajectories for the empirical basis functions, with a timestep of size 1. The mean square expectation of the numerical solution was computed at $x=0$ and compared to the true solution. Using only 4 or 5 empirical basis functions was inadequate for accuracy over the entire integration time, but 9 basis functions gave excellent agreement with the exact solution over the entire integration time, and adding additional basis functions does not appear to improve the accuracy. 

\begin{figure}[htb]
	\centering
	\begin{subfigure}[b]{0.48\textwidth}
    	\includegraphics[scale=0.35]{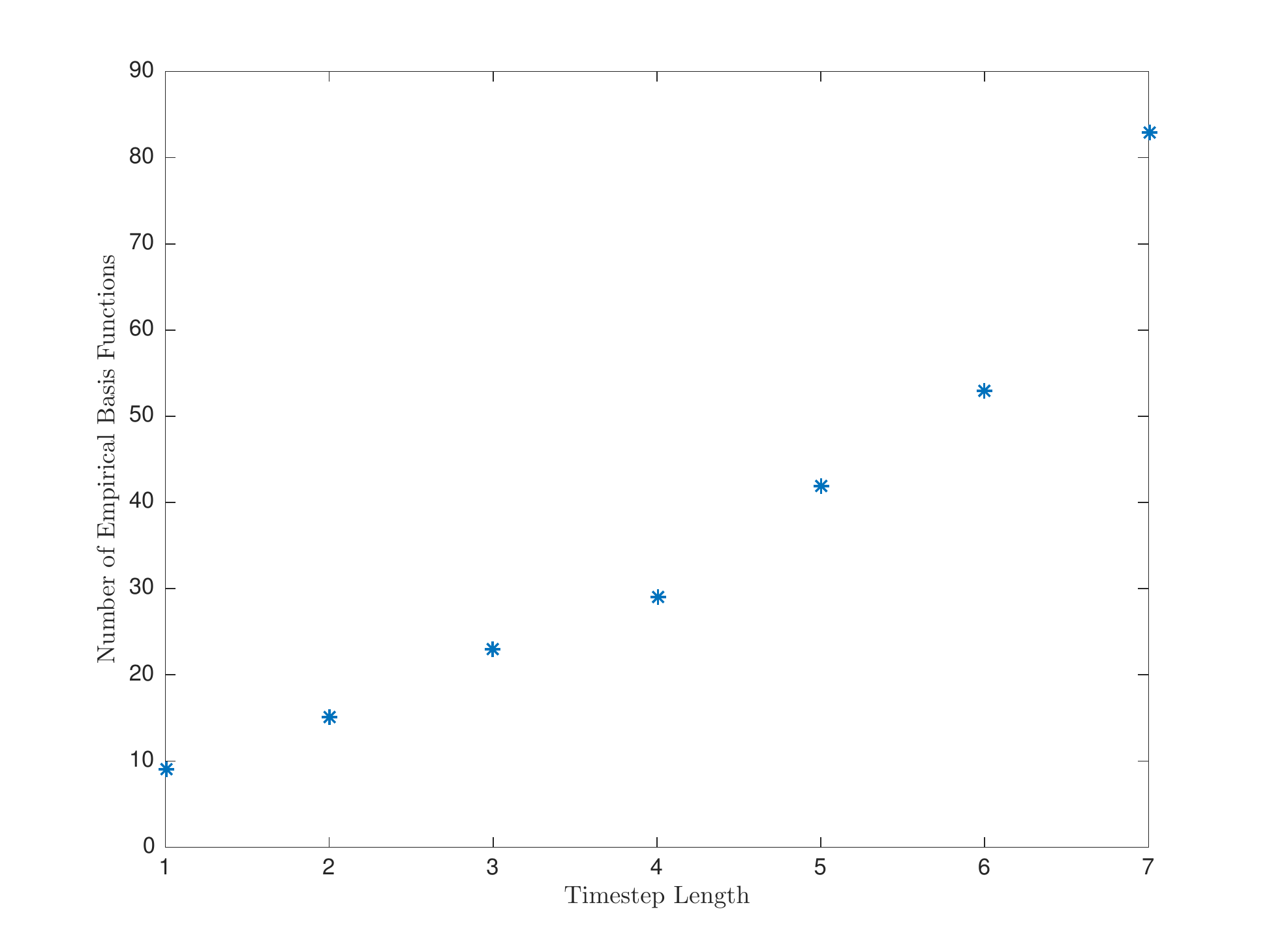}
    	\caption{\label{basis_size_vs_timestep}Number of empirical basis functions required for an accurate solution for varying timesteps. Note the fast growth rate.}
    \end{subfigure}\quad
	\begin{subfigure}[b]{0.48\textwidth}
  		\centering
    	\includegraphics[scale=0.35]{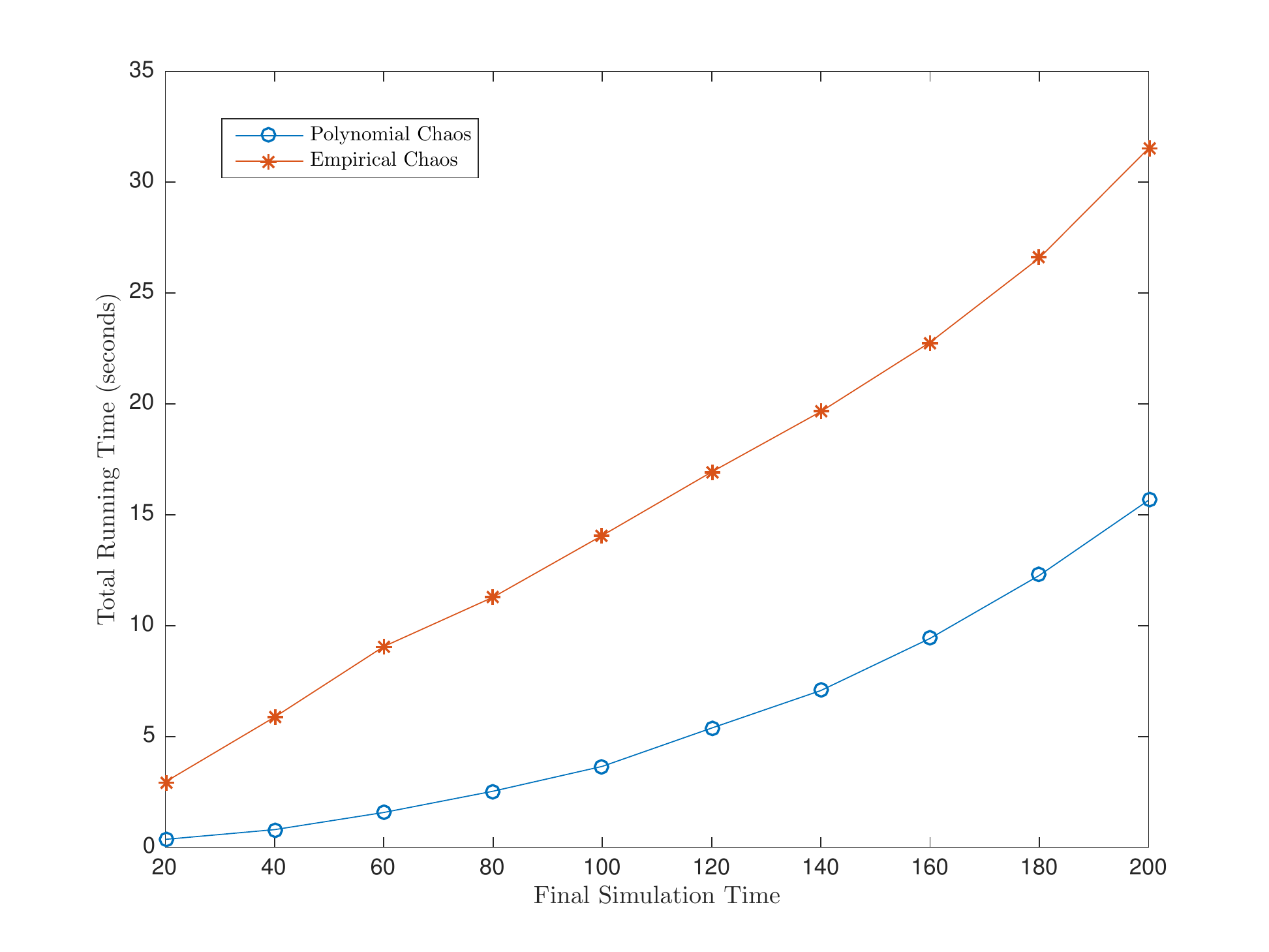}
    	\caption{\label{wave_1_comp}Comparison of total running times for empirical chaos expansion and gPC with stochastic Galerkin to the same accuracy.}
    \end{subfigure}
	\caption{Computational scaling of empirical chaos expansion with stochastic Galerkin as a function of timestep and simulation length for the one-dimensional wave equation \eqref{ow_wave}.}
\end{figure}

In general, choosing a larger timestep for the empirical chaos expansion results in growth of the required number of basis functions for each timestep. This is unsurprising, since we observe the same behavior with gPC; in order to accurately capture the solution behavior for long time periods, we need more and more basis functions. This presents an interesting optimization challenge for empirical chaos expansion, since a smaller basis results in a smaller stochastic Galerkin system, but a smaller timestep means that we must recompute the empirical bases more frequently. Recomputing the empirical basis involves sampling trajectories of the SPDE, computing the POD of the trajectories, constructing a new empirical basis, and projecting the data onto the new basis. Thus, a smaller timestep reduces the cost of solving the stochastic Galerkin system, but increases the cost of generating the empirical bases. For the model problem of this section, a timestep of size 4 appeared to be the optimal choice. Figure~\ref{basis_size_vs_timestep} shows the required number of basis functions for different choices of timesteps.

\subsection{Running Time}
As long as the number of empirical basis functions does not noticeably grow as the number of timesteps increases, we can expect the execution time for the empirical expansion method to scale linearly with the final integration time (see Figure~\ref{wave_1_comp}). When compared to the gPC method, we note the superlinear growth of the running time for gPC. However, for the timescales examined, gPC outperforms the empirical chaos method due to the fixed overhead cost of generating a new set of basis functions and performing the change of basis when applying an empirical chaos expansion.

Since the Legendre polynomials and the $A$ matrix in \eqref{ow_wave_pc} can be precomputed for gPC, the running time to construct them is not included in the gPC running time. However, as the order gets larger, the running time to construct the $A$ matrix is far larger than the actual time required to execute the gPC code, due to the $O(n^2)$ time required to compute all of the entries. In the simulations presented, the Legendre polynomials were generated symbolically and integrated exactly to avoid any error. In order to integrate out to time $200$ accurately, we needed $220$ basis functions, and the time required to compute the $A$ matrix was around 10 hours. In contrast, the empirical chaos method requires no precomputation and instead computes its basis functions directly from the sampled trajectories. 

\subsection{Basis Function Evolution}
\begin{figure}[htb]
	\centering
	\begin{subfigure}[b]{0.48\textwidth}
		\centering
		\includegraphics[scale=0.35]{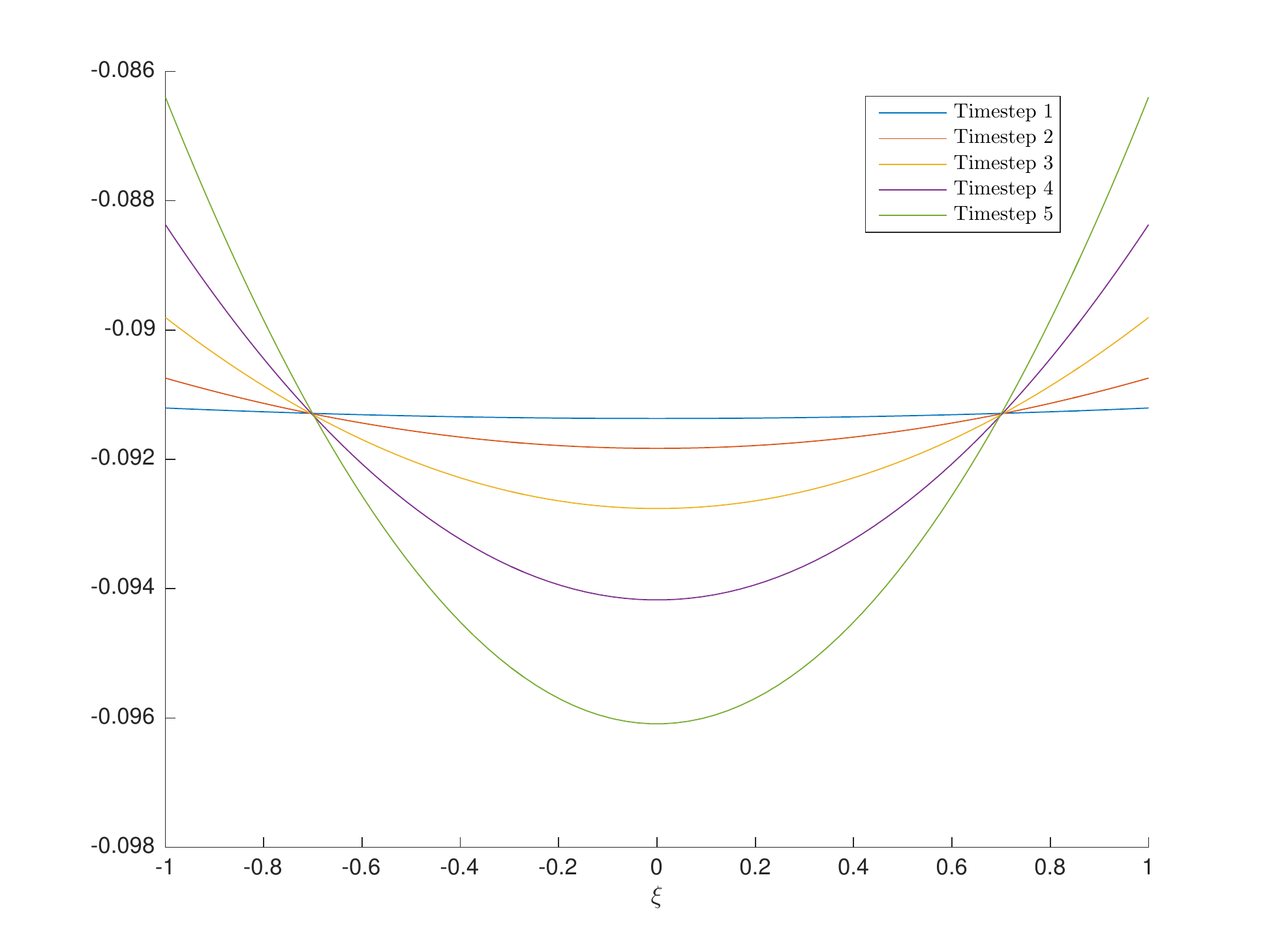}
		\caption{\label{basis_evo_1}Evolution of the first basis function from POD.}
	\end{subfigure}\quad
	\begin{subfigure}[b]{0.48\textwidth}
  		\centering
    	\includegraphics[scale=0.35]{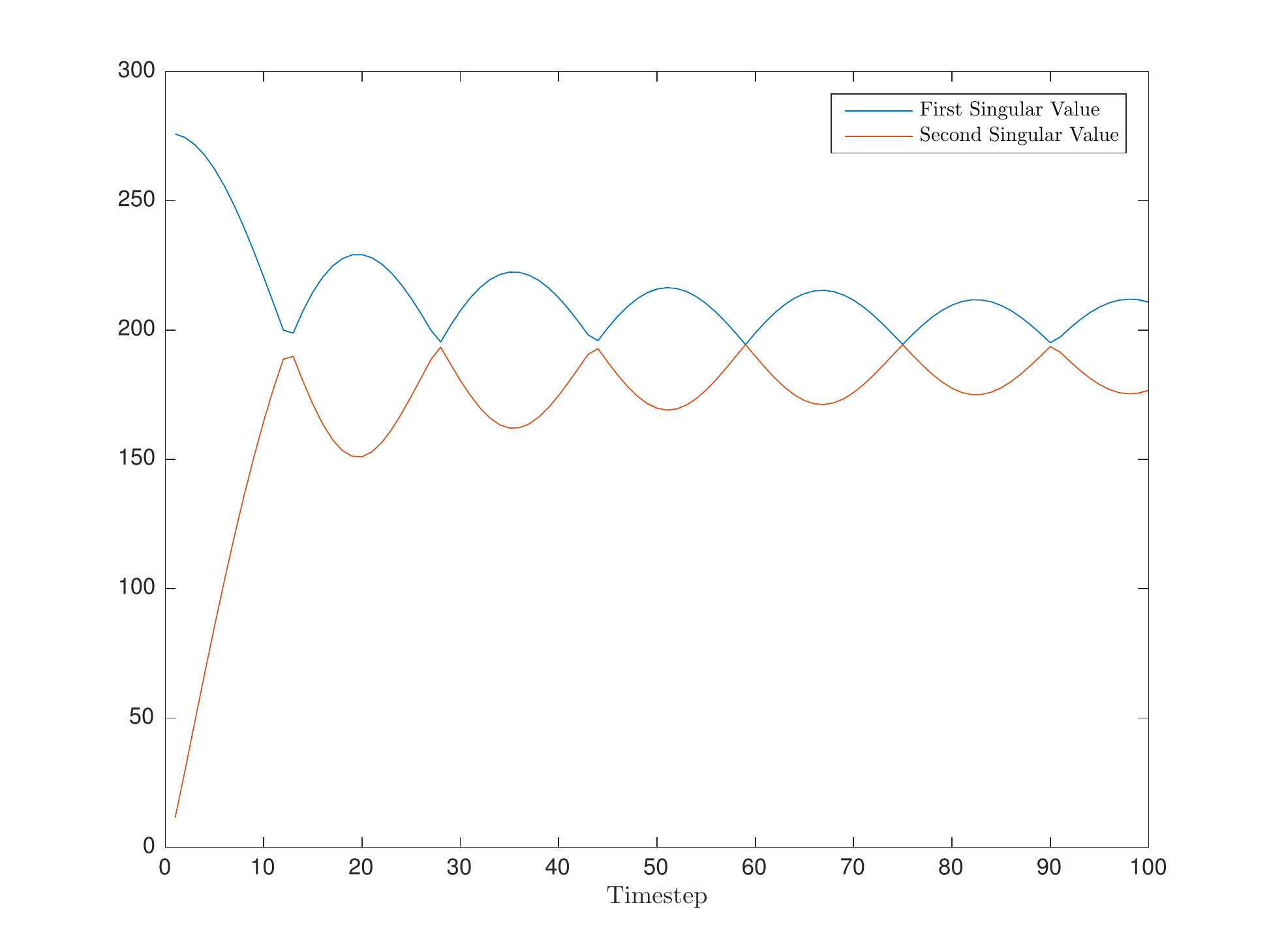}
    	\caption{\label{ow_wave_svs}Evolution of the magnitude of the first two singular values from POD.}
    \end{subfigure}
    \caption{Evolution of the basis functions and singular values from POD in the solution to \eqref{ow_wave} using empirical chaos expansion.}
\end{figure}

Since a new set of basis functions is computed at every timestep in the empirical expansion, it is natural to examine how they change over time. In order to closely monitor their evolution, we set the timestep for the wave equation to $0.1$, and examine the values of the basis functions over time. Figure~\ref{basis_evo_1} shows the basis function that corresponds to the largest singular value from POD for the first 5 timesteps (each of size 0.1). 

\begin{figure}[htb]
  \centering
    \includegraphics[scale=0.5]{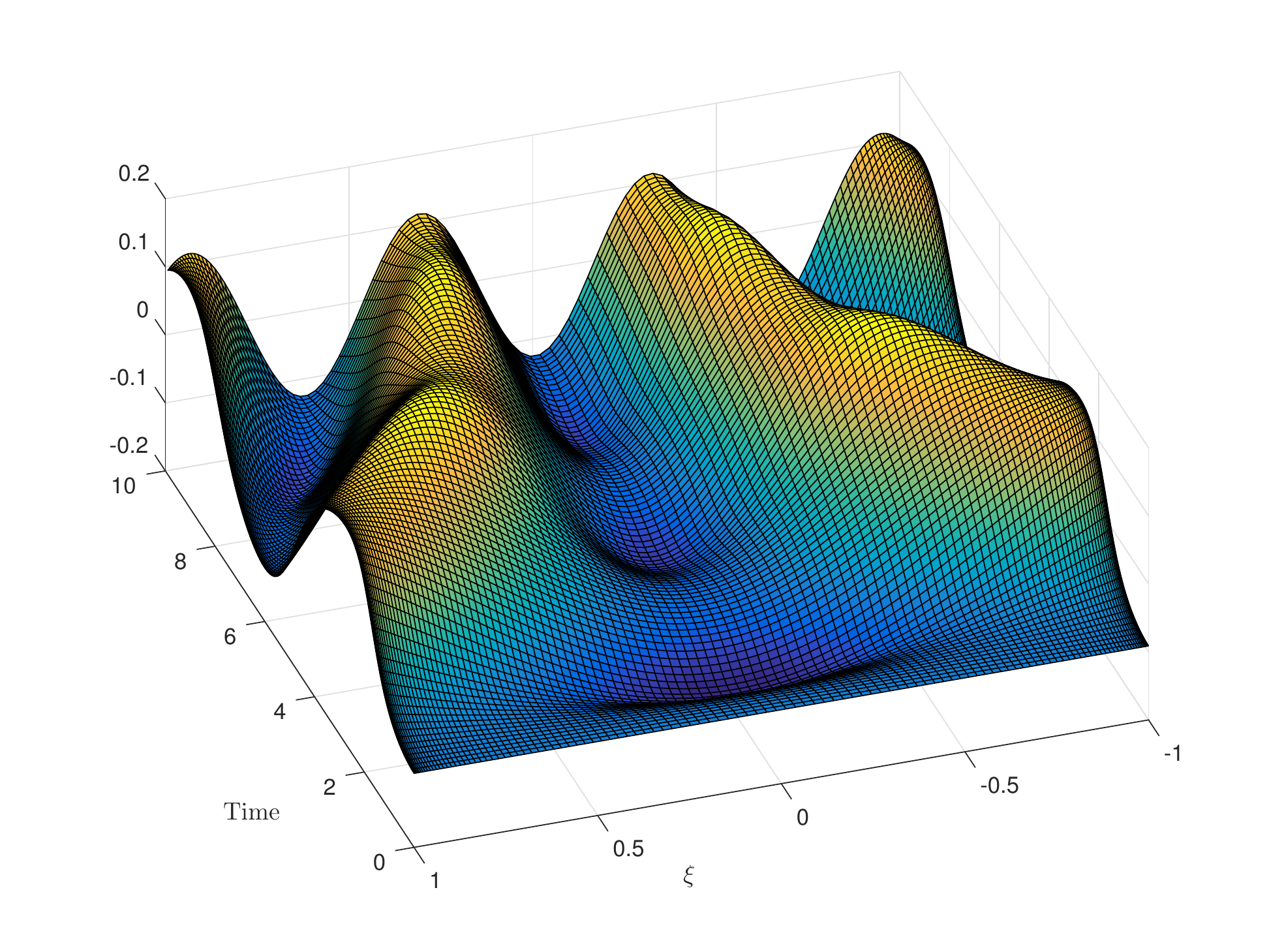}
    \caption[One-dimensional wave basis function evolution]{Evolution of the first basis function from POD for the first 100 timesteps in the solution to \eqref{ow_wave} using empirical chaos expansion.}
    \label{basis_evo_total}
\end{figure}

The basis function evolves smoothly over time, but at singular value crossings it becomes associated with the second largest singular value. Figure~\ref{ow_wave_svs} shows the magnitudes of the first and second singular values from POD for 100 timesteps, and there are multiple crossings. Figure~\ref{basis_evo_total} gives a three-dimensional view of the smoothness of the basis function evolution. Later, we will explore how to take advantage of this smooth evolution.

\section{Advection-Reaction Equation\label{adv_reac_sec}}
The gPC method outperforms the empirical chaos expansion even though it uses a far larger polynomial basis for the wave equation. There are two reasons for this. First, the stochastic Galerkin PDE system was relatively straightforward to solve, even in high-dimensions. Additionally, the polynomial chaos expansions rely on precomputing the $A$ matrix in \eqref{ow_wave_pc} since it always uses the fixed set of Legendre polynomials as its orthogonal basis. This accounts for much of the efficiency observed in the previous section.

Consider the advection-reaction equation,
\[
u_t(x,t,\xi) = \xi u_x(x,t,\xi) + f(x,u).
\]
If $f(x,u)$ is nonlinear and not a simple polynomial such as $u^2$, the stochastic Galerkin system is more difficult to evaluate than the original system. If the basis is large, then the stochastic Galerkin system for standard polynomial chaos becomes computationally intractable, and the benefit of empirical chaos expansion manifests itself. By keeping the basis small, the stochastic Galerkin system remains relatively inexpensive even when the final integration time is large.

In the subsequent sections, we consider the model advection-reaction SPDE,
\begin{align}
u_t(x,t,\xi) &= \xi u_{x}(x,t,\xi) + 0.1\abs{u}^{\frac{1}{2}},\quad 0\le x \le 2\pi, \quad t\ge0, \label{adv_reac_eq}\\
u(x,0,\xi) &= \cos(x) + \frac{3}{2}\label{adv_reac_bc},
\end{align}
with periodic boundary conditions, and $\xi$ is uniform on $[-1, 1]$. Note that
\begin{align*}
\iprod{f,g} \coloneqq \int_{-1}^1 f(\xi)g(\xi)(1/2)d\xi, \qquad \E{f} \coloneqq \int_{-1}^1 f(\xi)(1/2)d\xi.
\end{align*}
We do not have an analytic solution for this SPDE, so we run Monte Carlo simulations until they converge in order to obtain a reference solution to check the accuracy of our methods.

\subsection{Polynomial Chaos Expansion}
When using gPC, we choose the normalized Legendre polynomials as our basis functions. Let $L^i$ be the $i$-th normalized Legendre polynomial and apply stochastic Galerkin as in Section~\ref{section_pc}. The details are similar to Section~\ref{wave_gpc}. Letting
\begin{align*}
A_{ji} = \E{\xi L_jL_i},\quad
\hat{f}_j = 0.1\E{\abs{\sum\nolimits_{i=1}^\infty \hat{u}^i_x(x,t)L^i(\xi)}^{\frac{1}{2}} L^j(\xi)},\quad
\hat{u} =
\begin{pmatrix}
    \hat{u}^1(x,t)\\
    \hat{u}^2(x,t)\\
    \vdots\\
\end{pmatrix},
\end{align*}
the gPC system can be written as
\begin{equation}
\hat{u}_t = A\hat{u}_x + \hat{f}\label{adv_reac_sys}.
\end{equation}
Since the initial condition is deterministic,
\[
\hat{u}^1(x,0) = \cos(x)+ \frac{3}{2},\ \mathrm{and}\ \hat{u}^i(x,0)=0\ \mathrm{for}\ i>1.
\]
We can truncate the infinite system to finite $N$ and solve the deterministic system of PDEs. 

The stochastic Galerkin system in \eqref{adv_reac_sys} is more difficult to solve than the original SPDE system in \eqref{adv_reac_eq} due to the nonlinear $\hat{f}$ term. The expectations cannot be precomputed since they depend on $\hat{u}^i$ at each point in spacetime. We can compute $\hat{f}$ at each timestep with numerical quadrature, but it requires us to compute $N$ expectations, where $N$ is the number of basis functions. If we solve the stochastic Galerkin system using a method of lines discretization or some other finite-difference based approach, we will need to compute $N$ expectations for every timestep of the numerical solver. This becomes very expensive as $N$ grows. In order to compute the expectations in the $\hat{f}$ vector, we use a composite trapezoidal quadrature method that uses 300 Chebyshev nodes on the interval $[-1, 1]$. The Legendre polynomials and the $A$ matrix in \eqref{adv_reac_sys} are precomputed using symbolic arithmetic to avoid any error. 

\begin{figure}[htb]
	\centering
	\begin{subfigure}[b]{0.48\textwidth}
		\centering
    	\includegraphics[scale=0.35]{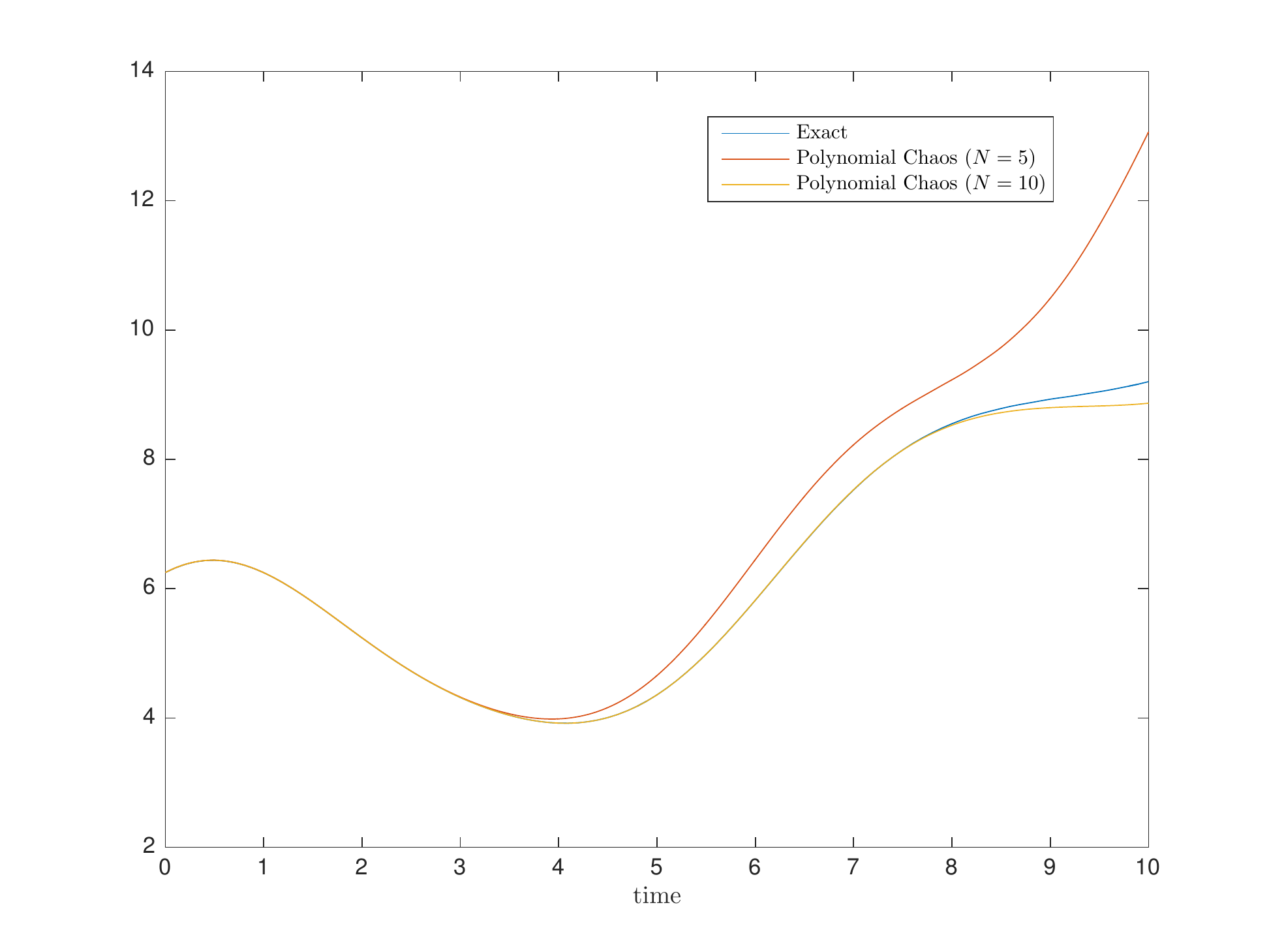}
    \caption{\label{adv_reac_tf10} Using $N$ polynomial basis functions.}
    \end{subfigure}\quad
    \begin{subfigure}[b]{0.48\textwidth}
		\centering
    	\includegraphics[scale=0.35]{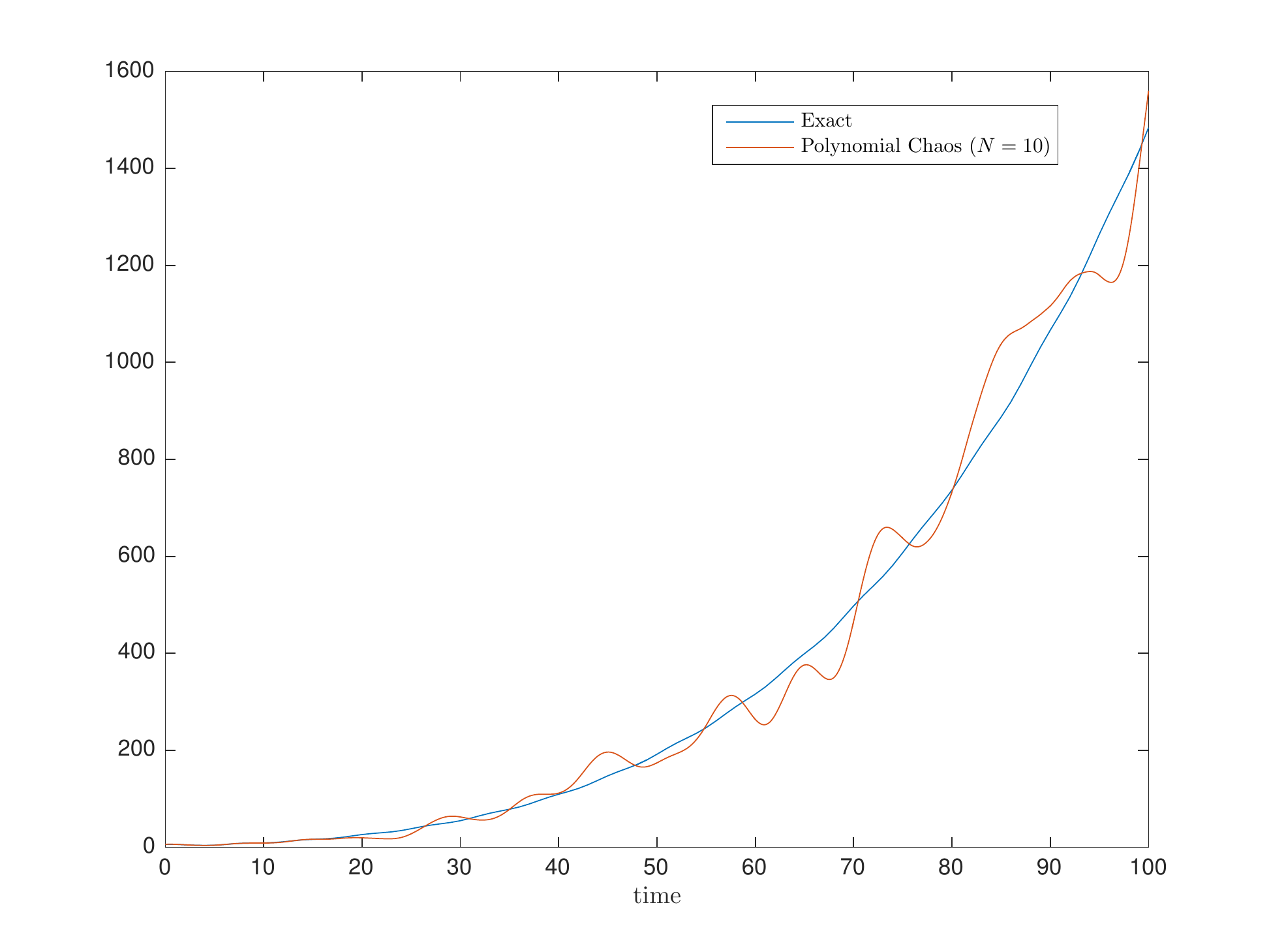}
    \caption{\label{adv_reac_tf100_nb10} Using Legendre polynomials up to order 10.}
    \end{subfigure}\\
    \begin{subfigure}[b]{0.48\textwidth}
		\centering
    	\includegraphics[scale=0.35]{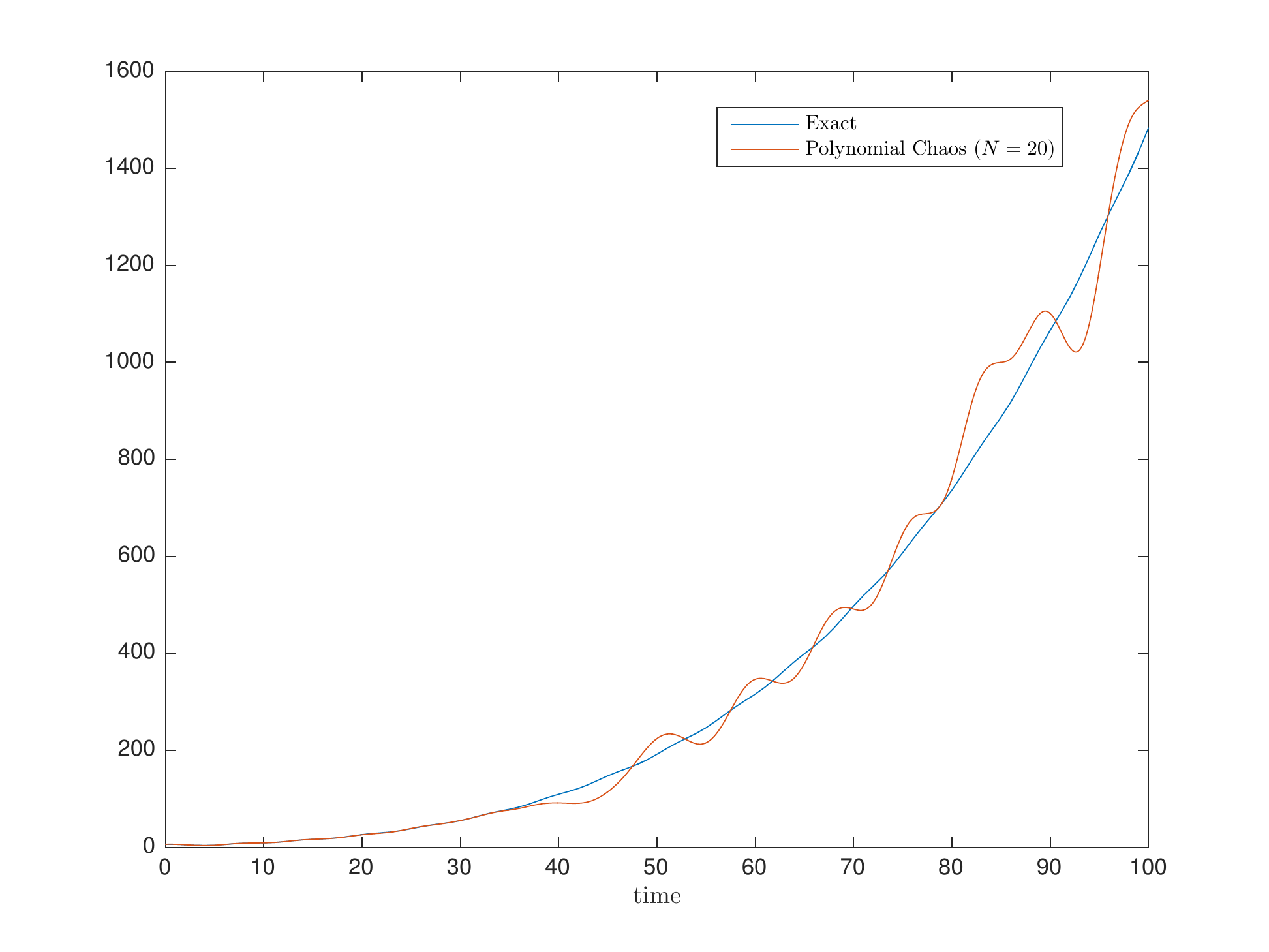}
    \caption{\label{adv_reac_tf100_nb20} Using Legendre polynomials up to order 20.}
    \end{subfigure}\quad
    \begin{subfigure}[b]{0.48\textwidth}
		\centering
    	\includegraphics[scale=0.35]{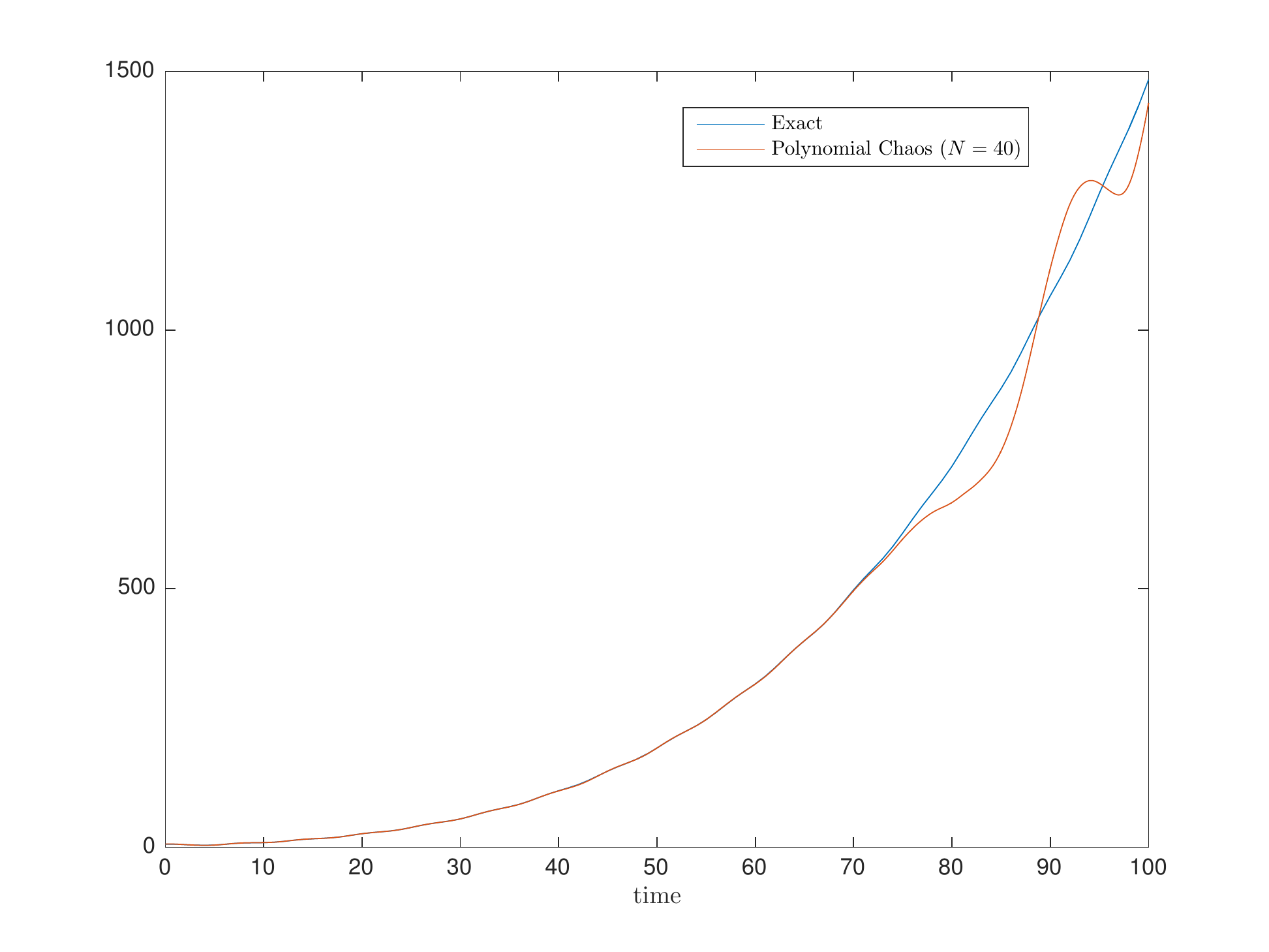}
    \caption{\label{adv_reac_tf100_nb40} Using Legendre polynomials up to order 40.}
    \end{subfigure}
    \caption{Mean square expectation $x=0$ of solution to \eqref{adv_reac_eq}, computed using gPC with stochastic Galerkin. The exact solution was obtained from 100,000 Monte Carlo realizations.}
\end{figure}

To examine the solution accuracy, we consider the mean square expectation at $x=0$, computed from the numerical solution of the stochastic Galerkin system in \eqref{adv_reac_sys}. Figures~\ref{adv_reac_tf10}, \ref{adv_reac_tf100_nb10}, \ref{adv_reac_tf100_nb20}, and \ref{adv_reac_tf100_nb40} show the numerical gPC results compared to the exact solution generated from 100,000 Monte Carlo simulations. For the short final time of $10$ in Figure~\ref{adv_reac_tf10}, a small basis performs well, but for the longer final time of $100$ in Figures~\ref{adv_reac_tf100_nb10}, \ref{adv_reac_tf100_nb20}, and \ref{adv_reac_tf100_nb40}, a much larger basis is required for accuracy. To numerically solve the differential equations, we use a method of lines discretization and a fourth-order Runge--Kutta integrator (\texttt{ode45} in MATLAB).

We include results up to polynomial order 40, but for higher orders, the stochastic Galerkin system became unstable. It is unclear if the instability was due to the numerical method used to solve the PDE system, or if the underlying system itself became ill-conditioned.

\subsection{Empirical Chaos Expansion}
If we let $\{\Psi^i(\xi)\}_{i=1}^N$ be the set of empirical basis functions, we can perform a stochastic Galerkin method as in Section~\ref{section_pc}. The details are similar to Section~\ref{wave_emp}. Letting
\begin{alignat*}{2}
    A_{ji} &= \E{\xi \Psi_j\Psi_i},&\qquad
    M_{ji} &= \E{\Psi_j \Psi_i},\\
    \hat{f}_j &= 0.1\E{\abs{\sum\nolimits_{i=1}^N \hat{u}^i_x(x,t)\Psi^i(\xi)}^{\frac{1}{2}} \Psi^j(\xi)},&\qquad
    \hat{u} &=
    \begin{pmatrix}
        \hat{u}^1(x,t)\\
        \hat{u}^2(x,t)\\
        \vdots\\
    \end{pmatrix},
\end{alignat*}
the empirical chaos system can be expressed as
\begin{equation}
    M\hat{u}_t = A\hat{u}_x + \hat{f}.
    \label{adv_reac_ec_sys}
\end{equation}
We follow the method in Section~\ref{emp_section} to construct an empirical basis over small time intervals by sampling trajectories and applying POD, and use these to compute $M$, $A$, and $\hat{f}$ in \eqref{adv_reac_ec_sys}. We use 300 Chebyshev nodes for $\xi$ in the interval $[-1,1]$. To fairly compare the computational cost with gPC, we compute the intermediate expectations of $\hat{f}$ using the same composite trapezoidal quadrature rule, and the same method of lines discretization and fourth-order Runge--Kutta integrator (\texttt{ode45} in MATLAB) to solve the stochastic Galerkin system.

\begin{figure}[htb]
	\centering
	\begin{subfigure}[b]{0.48\textwidth}
		\centering
    	\includegraphics[scale=0.35]{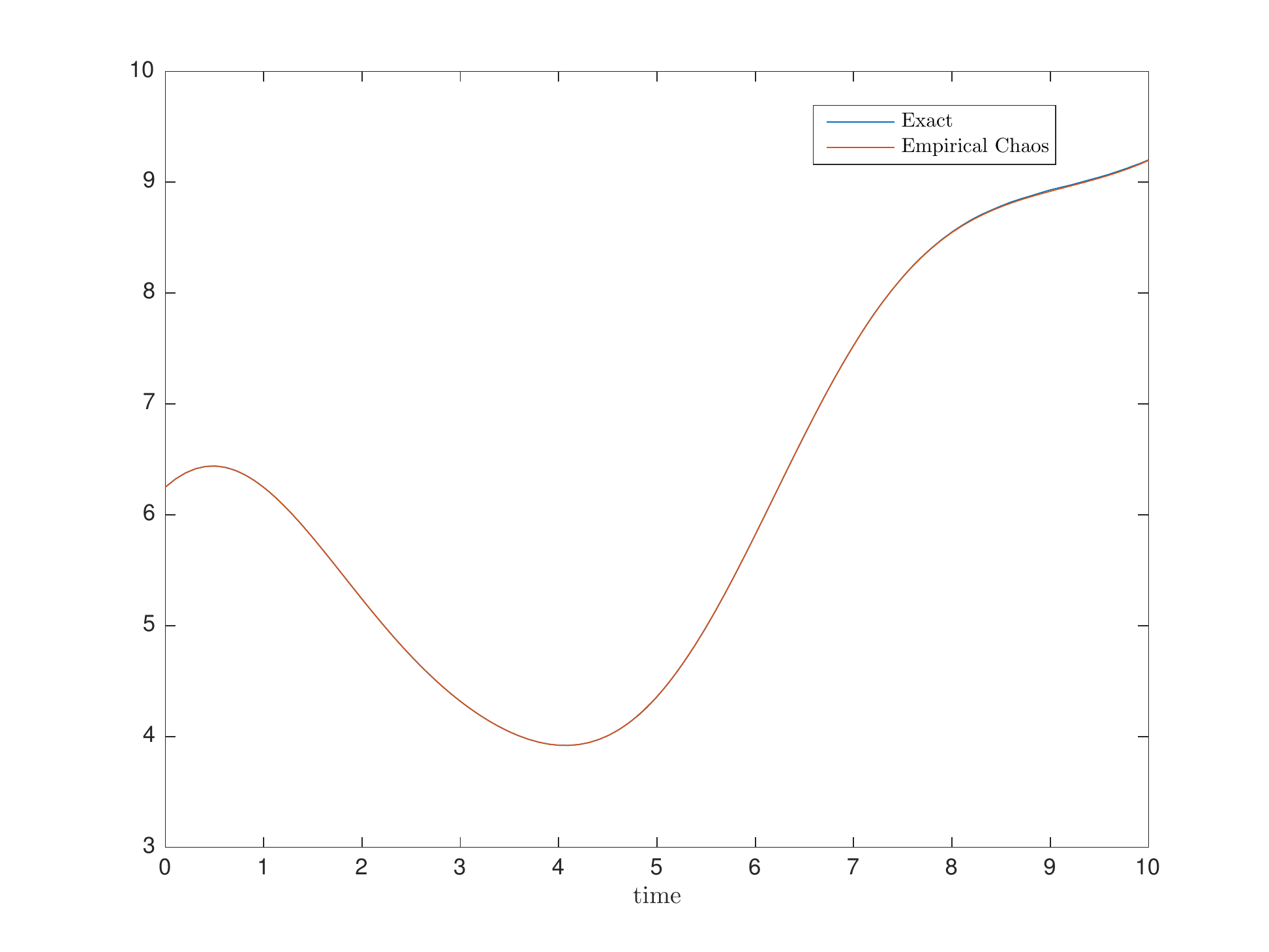}
    	\caption{\label{adv_reac_ec_tf10}$t\in[0,10]$.}
    \end{subfigure}\quad
    \begin{subfigure}[b]{0.48\textwidth}
		\centering
    	\includegraphics[scale=0.35]{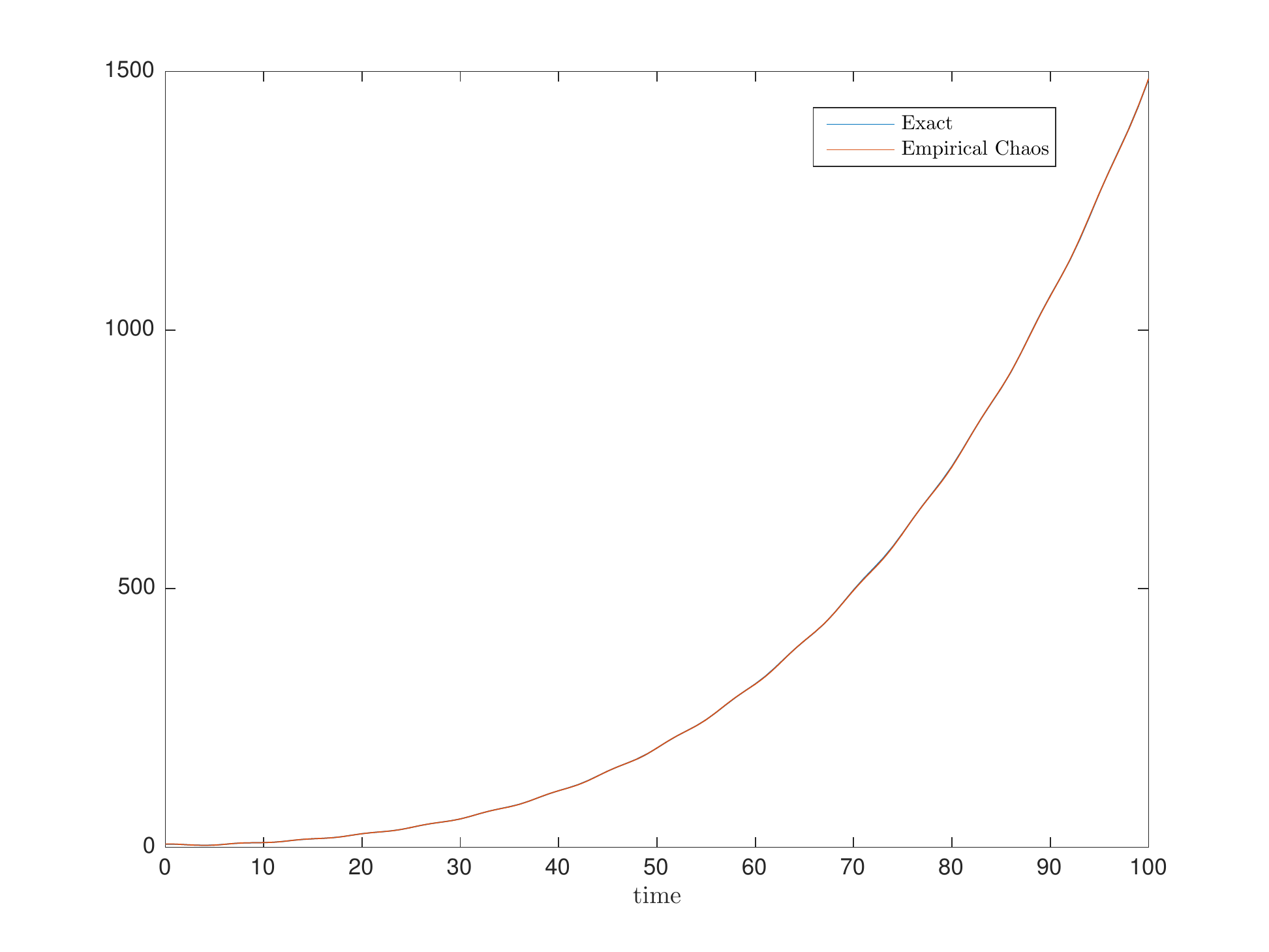}
    	\caption{\label{adv_reac_ec_tf100}$t\in[0,100]$.}
    \end{subfigure}
    \caption{Mean square expectation at $x=0$ of the solution to \eqref{adv_reac_eq} computed using an empirical chaos expansion with stochastic Galerkin. The timestep size was 2, the basis functions were computed by solving the deterministic equation at 300 Chebyshev nodes on the interval $[-1, 1]$, and the maximum number of basis functions used at a single timestep was 33.}   
\end{figure}

To examine the solution accuracy, we consider the mean square expectation at $x=0$, computed from the numerical solution to the stochastic Galerkin system in \eqref{adv_reac_ec_sys}. Figures~\ref{adv_reac_ec_tf10} and \ref{adv_reac_ec_tf100} compare empirical chaos expansion with the exact solution from 100,000 Monte Carlo simulations. The sampling points were 300 Chebyshev nodes in the interval $[-1,1]$.

\begin{figure}[htb]
	\centering
	\begin{subfigure}[b]{0.48\textwidth}
		\centering
		\includegraphics[scale=0.35]{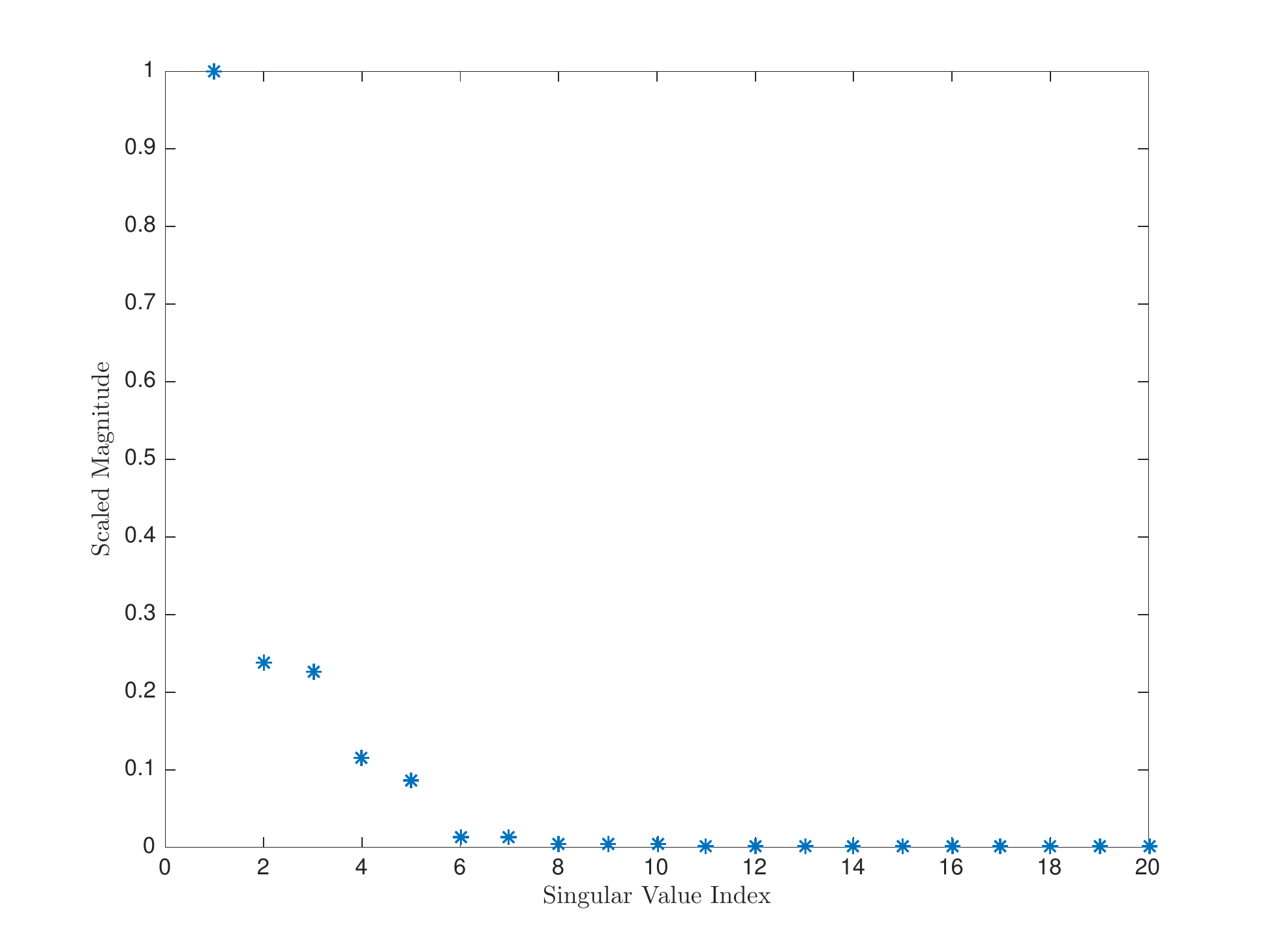}
    	\caption{\label{adv_reac_ec_svs6}Scaled singular values from POD for the third timestep.}
    \end{subfigure}\quad
	\begin{subfigure}[b]{0.48\textwidth}
  		\centering
    	\includegraphics[scale=0.35]{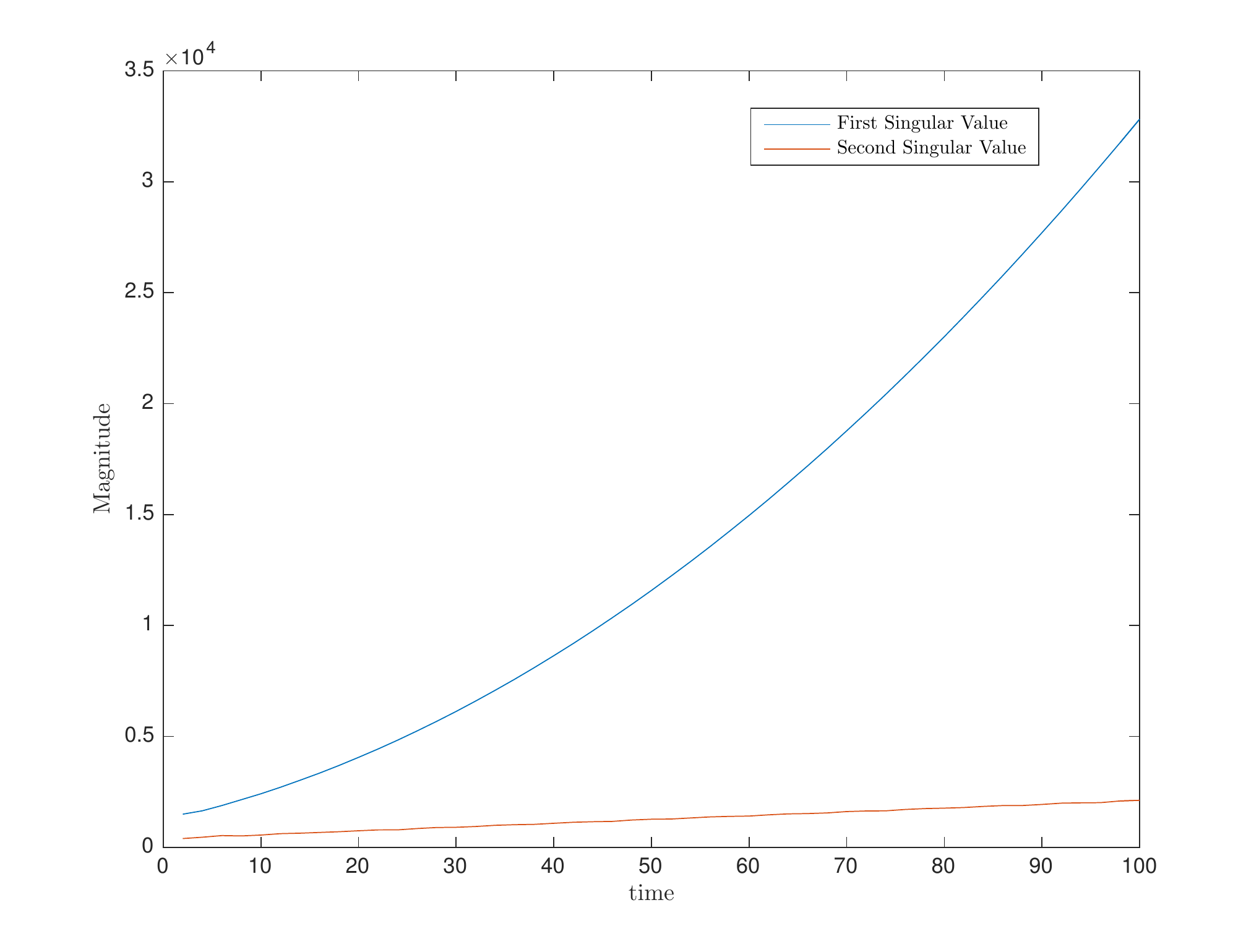}
    	\caption{\label{svs_ar_tf100}Evolution of the magnitude of the first two singular values from POD.}
    \end{subfigure}
	\caption{Singular values from POD in the solution to \eqref{adv_reac_eq} using empirical chaos expansion.}
\end{figure}

Just as before, we choose the number of empirical basis functions by examining the magnitude of the scaled singular values from POD. The overall number needed is higher for this SPDE than the previous ones, but there is still the same marked drop-off in the magnitude of the scaled singular values. Figure~\ref{adv_reac_ec_svs6} shows the magnitude of the singular values from the POD that was computed at $t=6$. There is a marked drop in the magnitude from the first to the second singular value, and by the sixth singular value the scaled magnitude is quite small. 

The initial mean square expectation fluctuates due to the presence of the advection term $\xi u_x$, but the reaction term $\abs{u}^{\frac{1}{2}}$ dominates the behavior out at long times. This can be seen when comparing the mean square expectations in Figures~\ref{adv_reac_ec_tf10} and \ref{adv_reac_ec_tf100}. In Figure~\ref{adv_reac_ec_tf10}, the solution is integrated to final time $10$, and there is a noticeable wave-like movement of the mean square expectation in time. In Figure~\ref{adv_reac_ec_tf100}, however, the solution is integrated to final time $100$, and though the wave-like motion is still present in the mean square expectation, it has been dominated by the total growth of the solution $u$ due to the reaction term. The singular values from POD also reflect this behavior. In Figure~\ref{svs_ar_tf100}, we plot the magnitudes of the first two singular values that were generated at each timestep of the empirical chaos expansion solver. As the time grows, the first singular value becomes completely dominant, which reflects the fact that the reaction term dominates the behavior at later times. 

\subsection{Running Time}
Since the number of empirical basis functions does not dramatically grow as the number of timesteps increases, the execution time for the empirical expansion method scales linearly with the final integration time (see Figure~\ref{wave_1_comp}). Note the superlinear growth of the gPC running time as a function of the final integration time. Although gPC outperforms the empirical chaos expansion for smaller integration times, it is outperformed by the empirical chaos expansion for final integration times larger than about 34.

We do not include the cost to precompute the $A$ matrix in \eqref{adv_reac_sys}, nor the cost to generate the Legendre polynomials for the polynomial chaos expansion. As mentioned previously, however, the entries of the $\hat{f}$ matrix cannot be precomputed, and the computational effort to compute these entries grows as $O(N^2)$, where $N$ is the number of basis functions. As the final simulation time increases, we must increase the number of polynomial basis functions in order to accurately solve the system, and this causes the computational cost to grow superlinearly.

\begin{figure}[htb]
  \centering
    \includegraphics[scale=0.35]{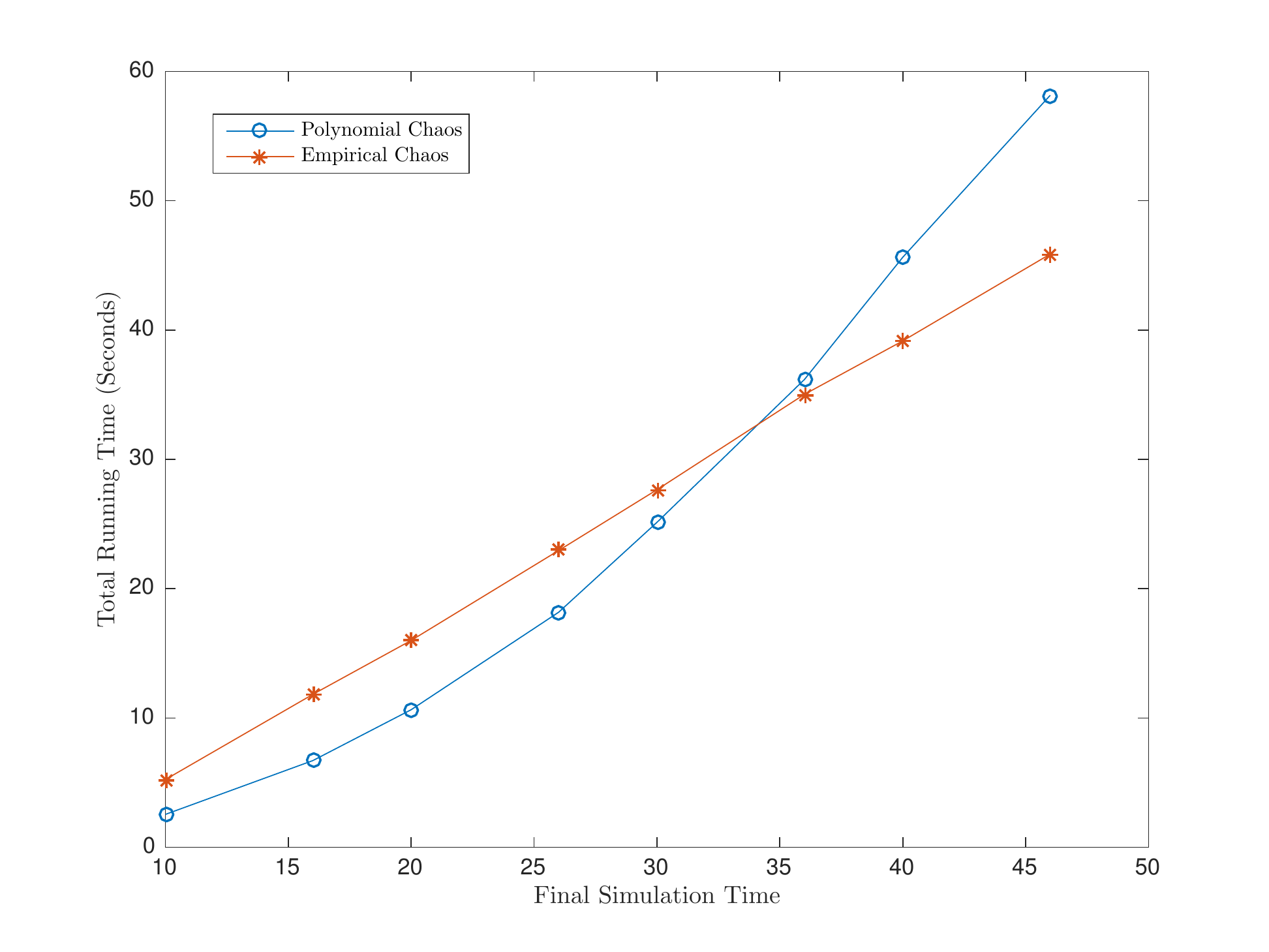}
    \caption[Advection-reaction running time comparison]{Comparison of total running times for solutions of \eqref{adv_reac_eq} computed to the same accuracy using empirical chaos expansion and gPC with stochastic Galerkin.}
    \label{avr_comp}
\end{figure}

\subsection{Basis Function Evolution}
\begin{figure}[htb]
	\centering
	\begin{subfigure}[b]{0.48\textwidth}
		\centering
		\includegraphics[scale=0.35]{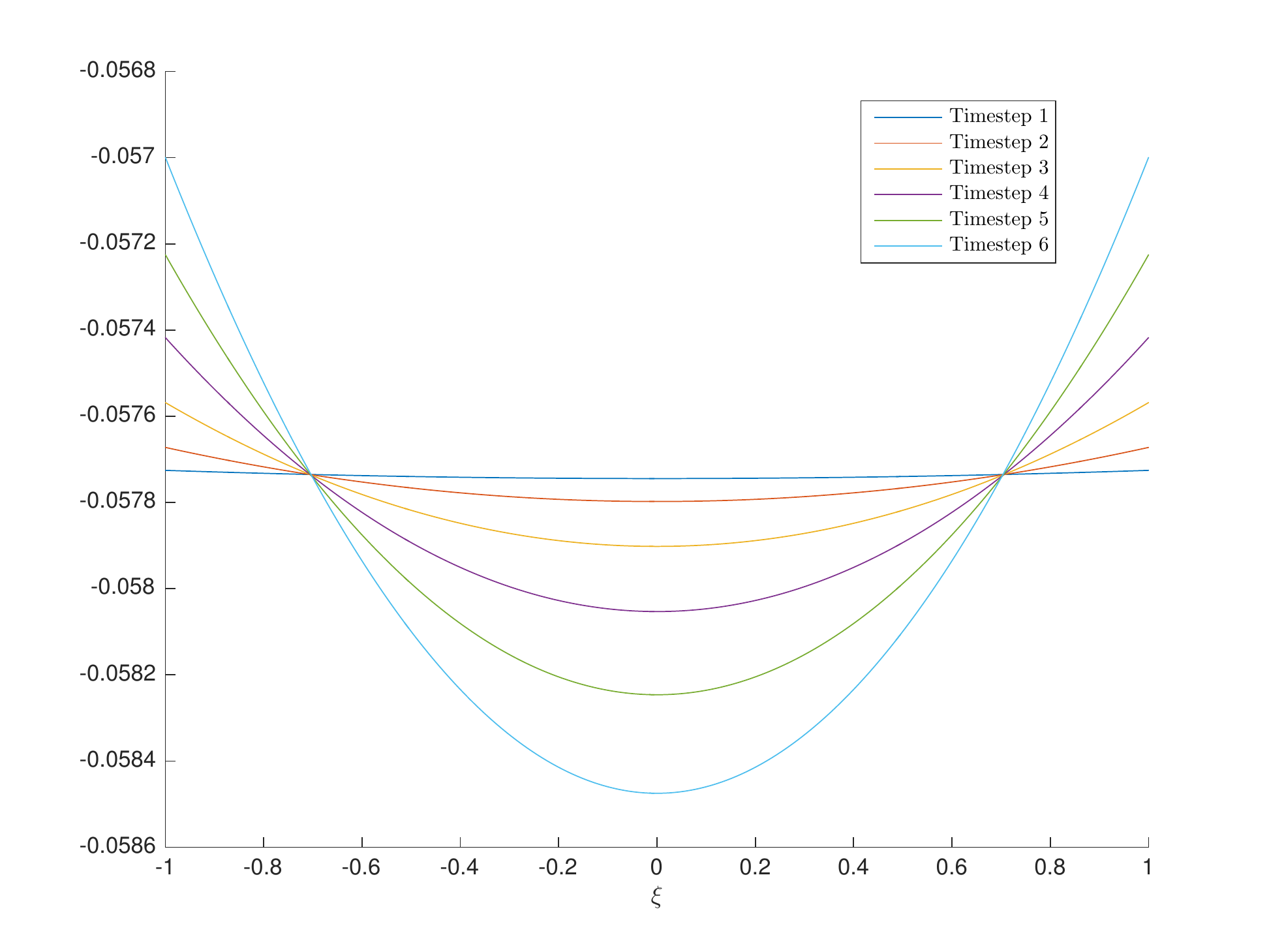}
    	\caption{\label{ar_basis_evo_sv1}Evolution of the first basis function from POD for the first five timesteps.}
	\end{subfigure}\quad
	\begin{subfigure}[b]{0.48\textwidth}
		\centering
		\includegraphics[scale=0.35]{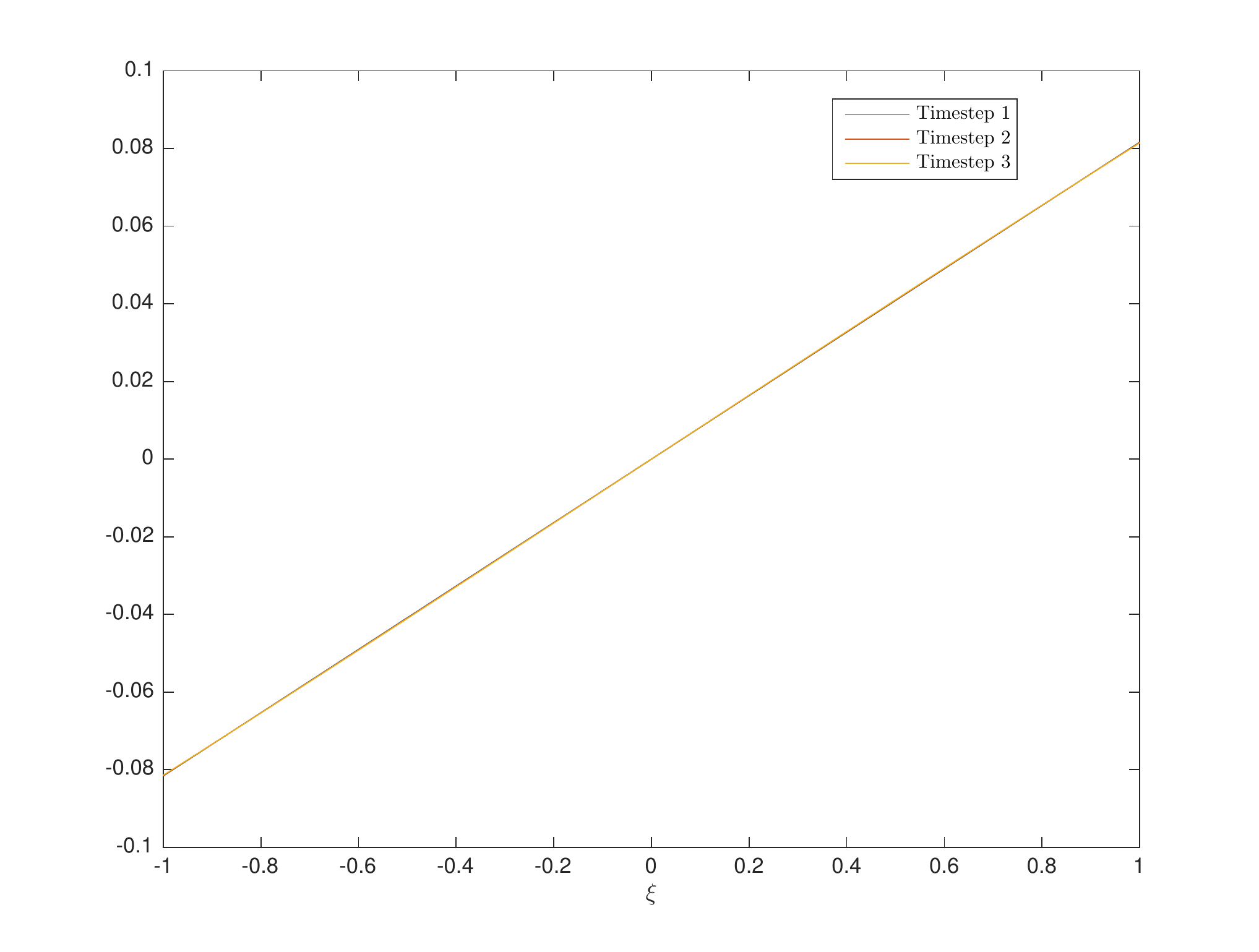}
    	\caption{\label{ar_basis_evo_sv2}Evolution of the second basis function from POD for the first three timesteps.}
	\end{subfigure}\\
	\begin{subfigure}[b]{0.48\textwidth}
		\centering
		\includegraphics[scale=0.35]{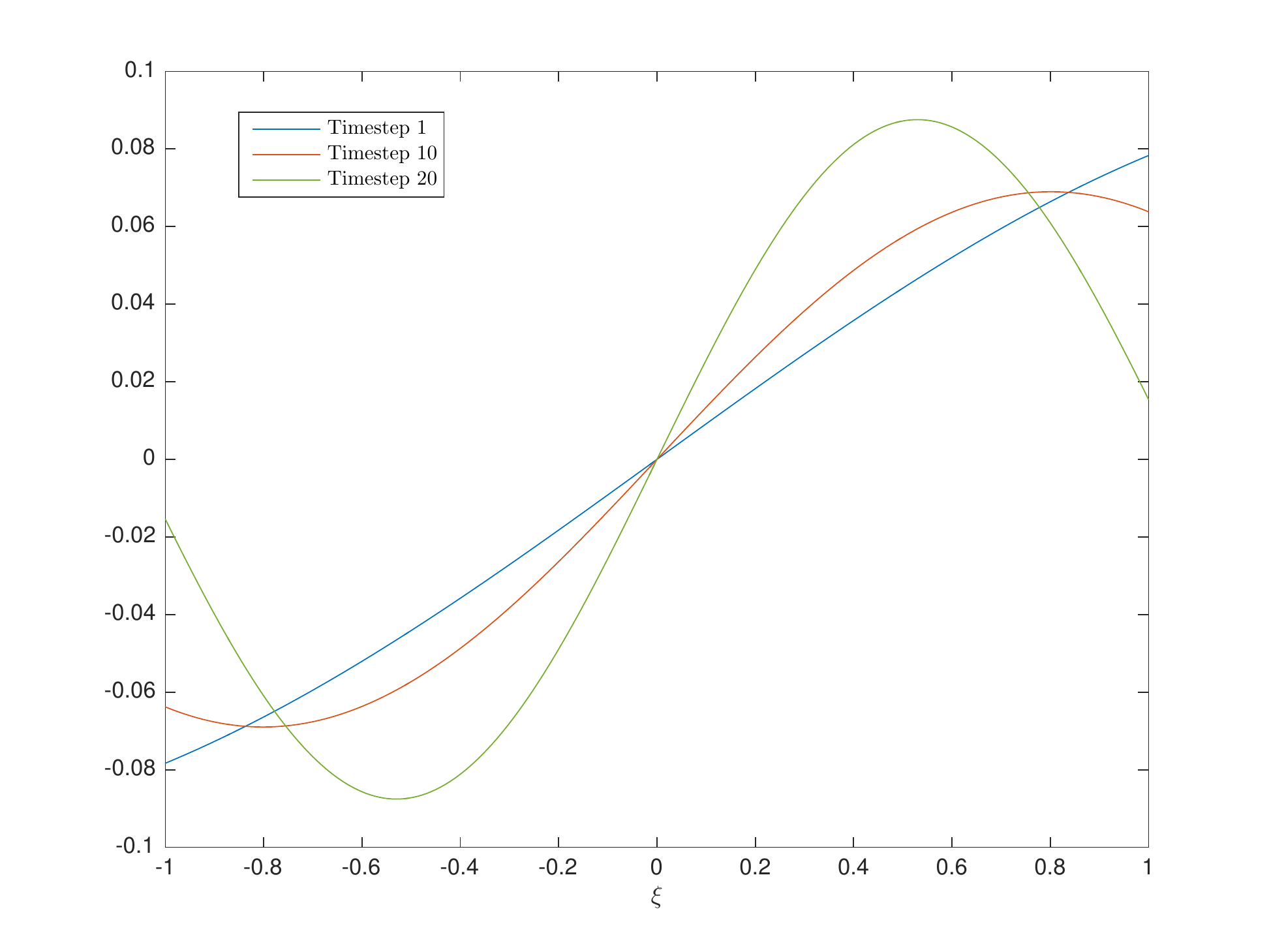}
    	\caption{\label{ar_basis_evo_sv2_ts10}Evolution of the second basis function from POD for the first 20 timesteps.}
	\end{subfigure}\quad
	\begin{subfigure}[b]{0.48\textwidth}
		\centering
		\includegraphics[scale=0.35]{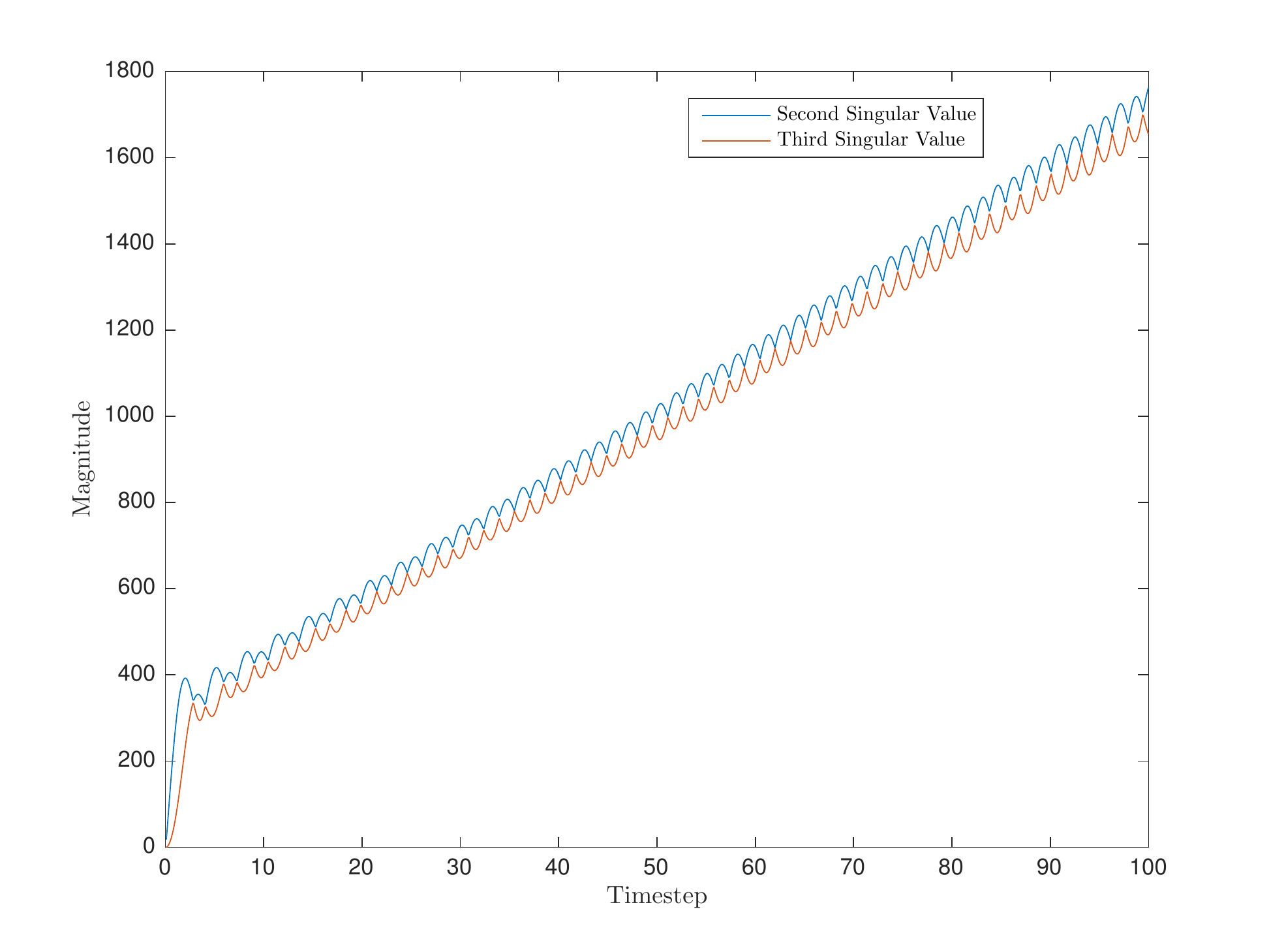}
    	\caption{\label{svs_ar_tf100_23}Evolution of the magnitude of the second and third singular values from POD.}
	\end{subfigure}
	\caption{Evolution of the basis functions from POD in the solution to \eqref{adv_reac_eq} using empirical chaos expansion.}
\end{figure}

In Figure~\ref{ar_basis_evo_sv1}, we examine the evolution of the basis function that corresponds to the largest singular value from POD for the advection-reaction equation~\eqref{adv_reac_eq} for the first $5$ timesteps (each of size $0.1$). The first basis function evolves smoothly over time, and since there are no crossings between the first and second singular value, it is straightforward to track its evolution. We also look at the second basis function evolution in Figures~\ref{ar_basis_evo_sv2} and \ref{ar_basis_evo_sv2_ts10}. From Figure~\ref{ar_basis_evo_sv2}, we see that the second basis function changes slowly, as is visually indistinguishable over the first 3 timesteps. The smooth evolution is more apparent in Figure~\ref{ar_basis_evo_sv2_ts10}, which tracks it over the first 20 timesteps. The second and third singular values have multiple crossings in Figure~\ref{svs_ar_tf100_23}, which means that, depending on the time, the second basis function is associated with either the second or third singular value.

\section{Empirical Basis Evolution}
The smooth evolution of the empirical basis functions derived from POD suggests that the empirical basis functions can be computed through an analytic approach that involves the dynamics of the SPDE. In gPC methods, we construct a basis that is orthogonal to the distribution of the random variables, project the exact solution onto this basis, and solve for the deterministic coefficients which are functions of spacetime. To examine how the dynamics of a SPDE influence the empirical basis functions for the random space, we consider an expansion in terms of a time-dependent basis for the random space,
\begin{equation}
u(x,t,\xi) \approx \sum\nolimits_{i=1}^N \hat{u}^i(x) \Psi^i(\xi, t), \label{basis_evo_eq}
\end{equation}
and match this with the empirical chaos expansion~\eqref{empirical_chaos_expansion} at time $t_0$, i.e., $\Psi^i(\xi,t_0)=\Psi^i(\xi)$ and $\hat{u}^i(x)=\hat{u}^i(x,t_0)$. The spatial coefficients are fixed and the random basis is time-dependent.

In order to evolve the basis functions in time, we consider a Galerkin method, this time integrating over the spatial component. We multiply by a test function $\hat{u}^j$, then integrate over the spatial domain, $R$. For the general model problem in \eqref{pc_model}, we obtain
\begin{align}
\sum\nolimits_{i=1}^N \int_{R}\hat{u}^i(x) \hat{u}^j(x)\ dx \Psi^i(\xi, t) &= \int_{R} L\left(\sum\nolimits_{i=1}^N \hat{u}^i(x) \Psi^i(\xi, t), x, t, \xi \right) \hat{u}^j(x)\ dx.\label{evo_sys}
\end{align}
In practice, this system can usually be simplified, but it depends on the form of the differential operator $L$. This provides a DAE that describes the evolution of the basis functions of the random variable in time. The general approach we follow is to generate a set of empirical basis functions over a short time interval using POD, solve the resulting stochastic Galerkin system, and use the final values of the $\hat{u}^i(x,t)$ coefficients as the fixed spatial basis in~\eqref{evo_sys}. We use the empirical basis functions as the initial values of the $\Psi^i(\xi, t)$ coefficients, and then solve~\eqref{evo_sys} over a short time interval and use the final values of the $\Psi^i(\xi,t)$ coefficients as the new empirical basis for the next stochastic Galerkin step.

\subsection{One-Dimensional Wave Equation}
If we substitute \eqref{basis_evo_eq} into the one-dimensional wave equation~\eqref{ow_wave}, we get $\sum\nolimits_{i=1}^N \hat{u}^i(x) \Psi^i_t(\xi, t) = \xi\sum\nolimits_{i=1}^N \hat{u}^i_x(x) \Psi^i(\xi, t)$. Multiplying by a test function $\hat{u}^j$ and integrating over the spatial domain yields
\begin{align*}
	\sum\nolimits_{i=1}^N \int_0^{2\pi}\hat{u}^i(x)\hat{u}^j(x)\ dx \Psi^i_t(\xi, t) &= \xi\sum\nolimits_{i=1}^N \int_0^{2\pi}\hat{u}^i_x(x)\hat{u}^j(x)\ dx \Psi^i(\xi, t).
\end{align*}
Letting
\[
\Psi = \begin{pmatrix}
\Psi^1\\
\Psi^2\\
\vdots\\
\Psi^N
\end{pmatrix},
\quad
\hat{u} = \begin{pmatrix}
\hat{u}^1\\
\hat{u}^2\\
\vdots\\
\hat{u}^N
\end{pmatrix},
\quad
A_{ji} = \int_0^{2\pi} \hat{u}^i(x)\hat{u}^j(x)\ dx,
\quad
M_{ji} = \int_0^{2\pi} \hat{u}^i_x(x)\hat{u}^j(x)\ dx,
\]
implies that
\begin{equation}
A\Psi_t(\xi, t) = \xi M \Psi(\xi, t),\label{owwave_basis_evo_eq}
\end{equation}
which is a DAE that describes the evolution of basis functions in time. The general approach is to first generate a set of empirical basis functions up to a particular time $t^*$, using the method of Section~\ref{emp_section}. Then, the set of empirical basis functions is evolved in time by solving \eqref{owwave_basis_evo_eq}. In particular, if the matrix $A$ is invertible, then the exact solution to \eqref{owwave_basis_evo_eq} is
\begin{equation}
    \Psi(\xi, t) = \exp{\left[\xi A^{-1} M(t-t_0)\right]}\Psi(\xi, t_0), \label{owwave_basis_exp_op_eq}
\end{equation}
where $\exp$ is the matrix exponential operator. In practice, $A$ is singular. To deal with this, we can decompose the full system into subsystems where the $A$ matrix is nonsingular. To accomplish this, we choose the largest nonsingular square submatrix of $A$ whose upper left corner is the $(1,1)$ entry of $A$. If that submatrix has lower right corner $(k,k)$ where $k<n$, then we repeat this procedure by choosing another square matrix whose upper left corner is $(k+1, k+1)$, and so on. In practice, we calculate the condition number of the submatrices to determine if the matrix is close to singular. Then, each subsystem has a nonsingular $A$ matrix, and the matrix exponential in \eqref{owwave_basis_exp_op_eq} can be directly applied to update the basis functions.

\begin{figure}[htb]
	\centering
	\begin{subfigure}[b]{0.48\textwidth}
		\centering
    	\includegraphics[scale=0.35]{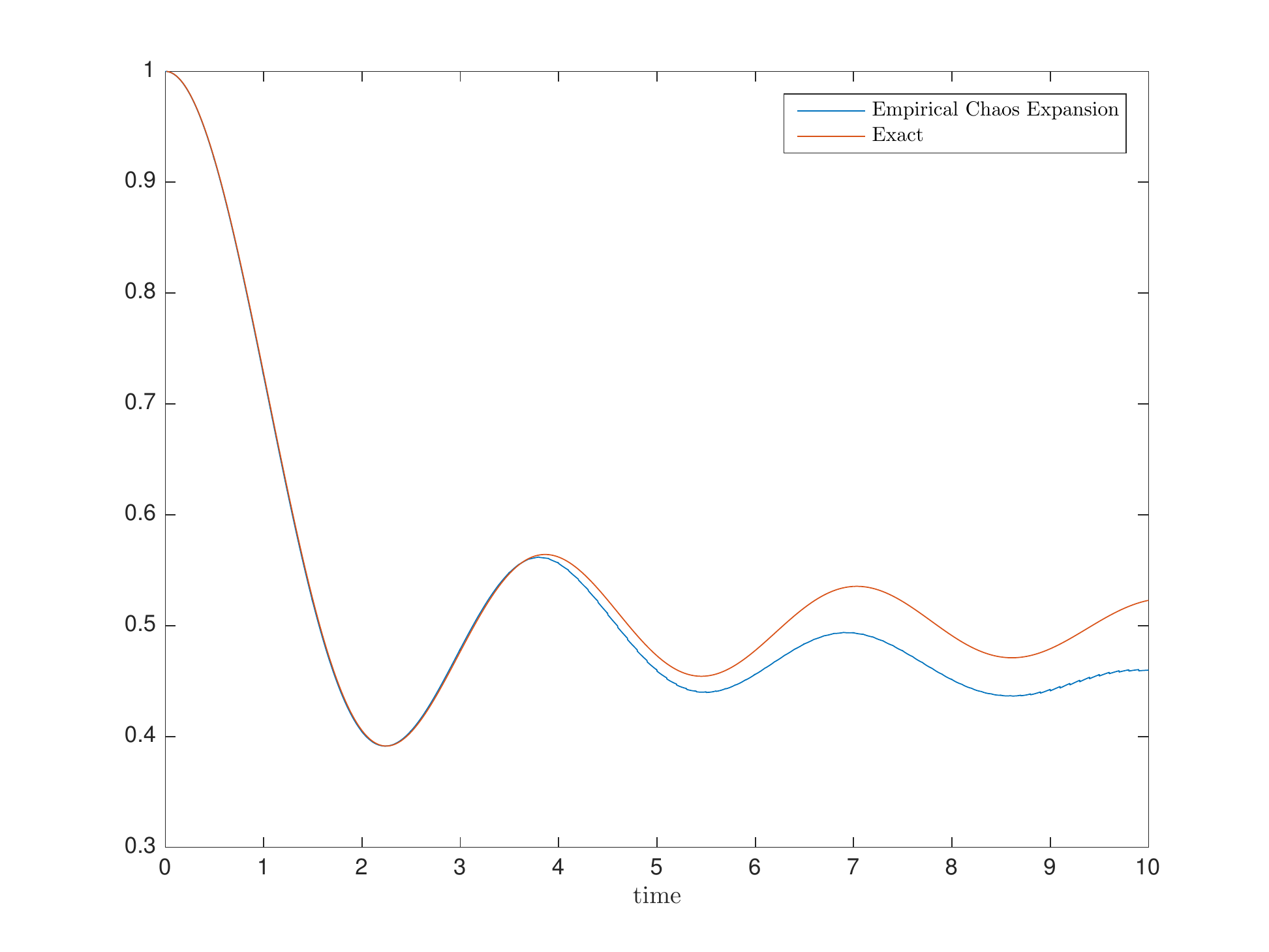}
    	\caption{\label{basis_evo_tf1_10} The basis functions up to time $t=1$ are empirically generated. For times $1 < t \le 10$, the empirical basis functions generated at $t=1$ are evolved using the exponential operator in \eqref{owwave_basis_exp_op_eq}.}
    \end{subfigure}\quad
    \begin{subfigure}[b]{0.48\textwidth}
		\centering
    	\includegraphics[scale=0.35]{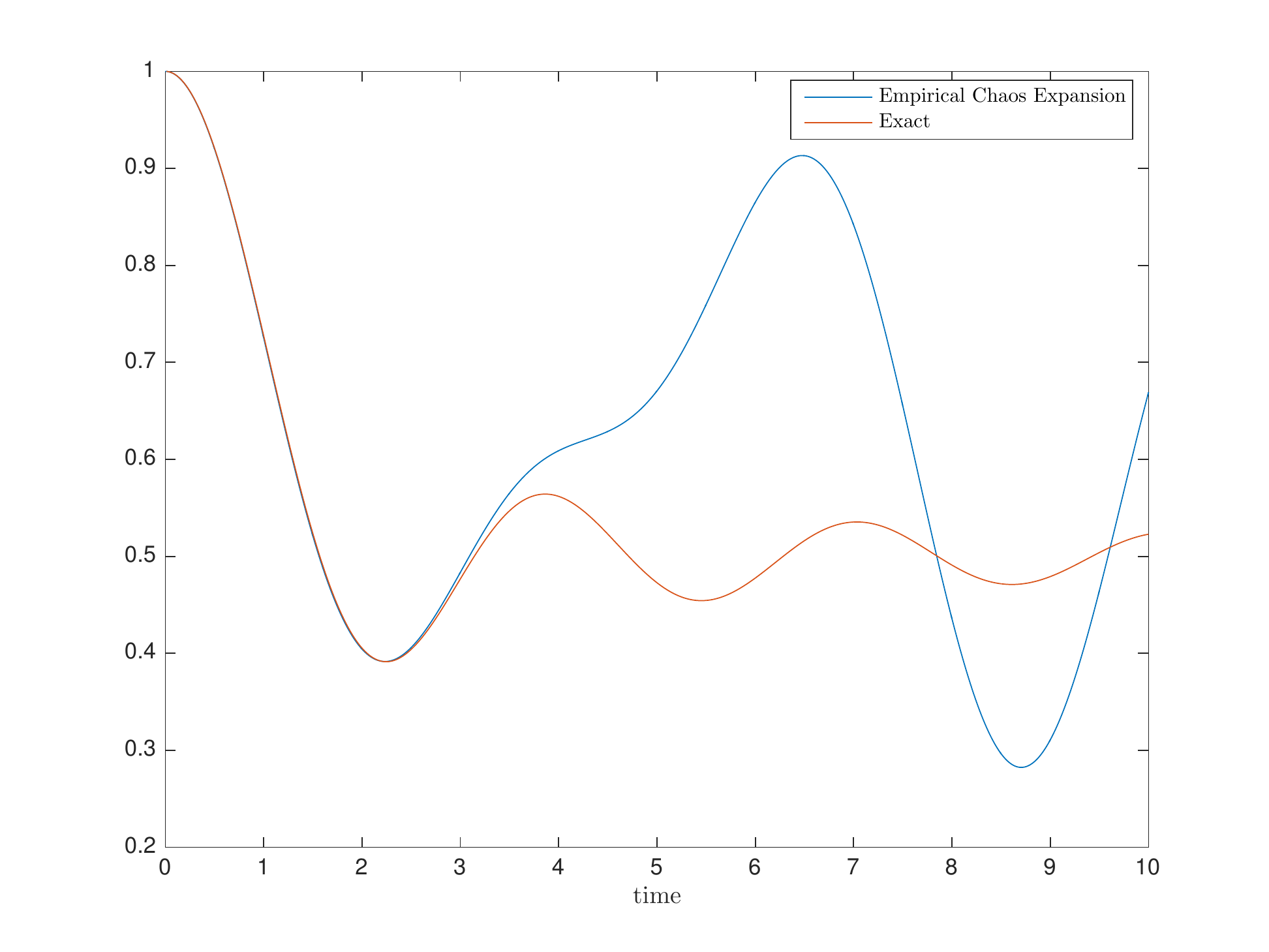}
    	\caption{\label{basis_fixed_tf1_10} The basis functions up to time $t=1$ are empirically generated. For times $1 < t \le 10$, the basis functions remain fixed at the empirical basis functions that were generated for time $t=1$.}
    \end{subfigure}\\
    \begin{subfigure}[b]{0.48\textwidth}
		\centering
    	\includegraphics[scale=0.35]{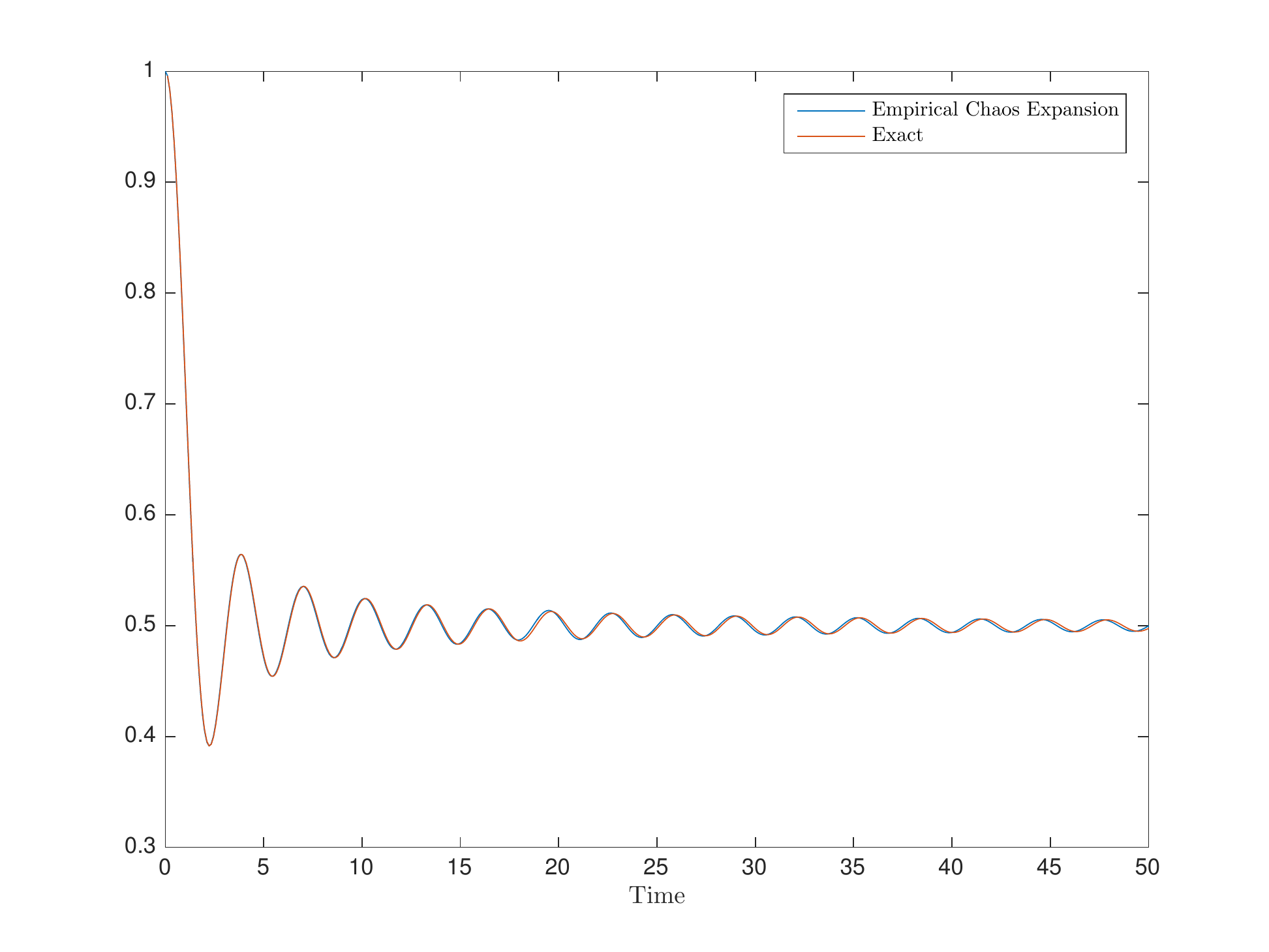}
    	\caption{\label{basis_evo_alt_tf50} The basis functions on even timesteps are empirically generated using POD. The basis functions on odd timesteps are evolved using the exponential operator in \eqref{owwave_basis_exp_op_eq}.}
    \end{subfigure}\quad
    \begin{subfigure}[b]{0.48\textwidth}
		\centering
    	\includegraphics[scale=0.35]{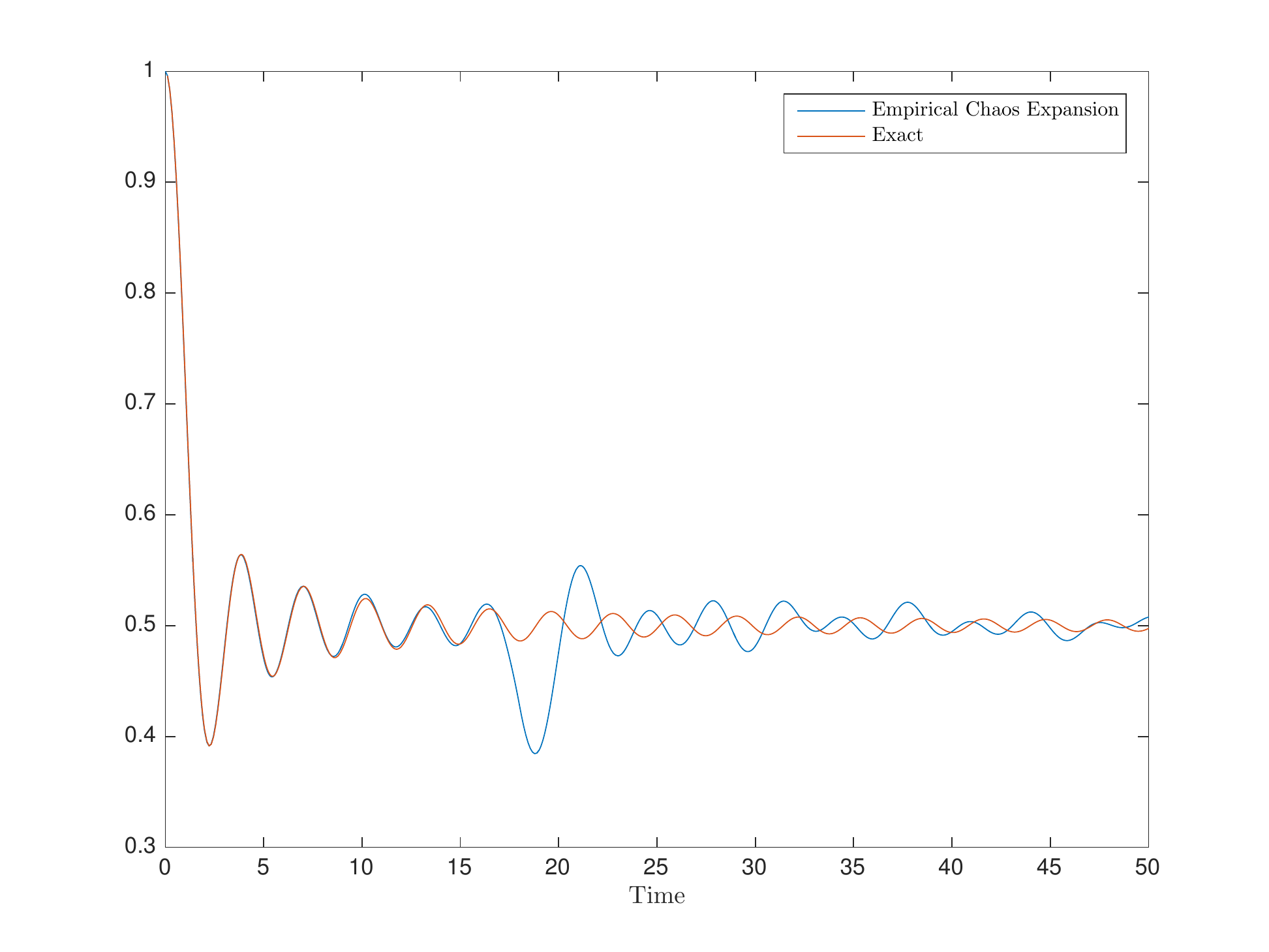}
    	\caption{\label{basis_noevo_alt_tf50} The basis functions on even timesteps are empirically generated using POD. The basis functions on odd timesteps are not evolved from their previous values.}
    \end{subfigure}
  \caption{Mean square expectation at $x=0$ of solution to \eqref{ow_wave}, computed using empirical chaos expansion with stochastic Galerkin.}
\end{figure}

Figure~\ref{basis_evo_tf1_10} illustrates this method applied to the one-dimensional wave equation~\eqref{ow_wave}. The basis functions up to time $t=1$ are empirically generated using the standard method from Section~\ref{owsection}. For times $1 < t \le 10$, the empirical basis functions that were generated at $t=1$ are evolved using the exponential operator from \eqref{owwave_basis_exp_op_eq}. The numerical solution closely matches the exact solution up to about $t=4$. Past that time, the numerical solution begins to diverge from the exact solution, but still follows the general pattern. For comparison, in Figure~\ref{basis_fixed_tf1_10}, the basis functions are empirically generated up to time $t=1$, and then the basis functions are kept fixed for times $1 < t \le 10$. When the basis functions are not evolved at all, the numerical solution matches the exact solution up to about $t=3$, but wildly diverges from the exact solution beyond that, in much the same way as the standard polynomial chaos expansion. 

The basis evolution operator involves a matrix exponential, and is potentially cheaper to compute than the trajectory sampling approach, which suggests the following general strategy:
\begin{enumerate}
    \item Compute sample trajectories and use POD to generate an empirical basis;
    \item Use the exponential basis evolution operator to update the basis functions for a few additional timesteps without resampling the trajectories;
    \item Repeat.
\end{enumerate} 

Although we do not resample the trajectories at each timestep, we need to use smaller timesteps when applying the matrix exponential operator. In the case of the one-dimensional wave, this results in a slight increase in the running time since the stochastic Galerkin system is already straightforward to solve. Figure~\ref{basis_evo_alt_tf50} illustrates the result of the alternating timestep strategy, where the first timestep generates empirical basis functions by sampling trajectories, the second timestep evolves the basis functions using the exponential operator from~\eqref{owwave_basis_exp_op_eq}, and so on. The solution retains the same accuracy as when we recomputed the basis at every timestep. In Figure~\ref{basis_noevo_alt_tf50}, we see that not using the exponential operator from~\eqref{owwave_basis_exp_op_eq} to update the basis functions on alternating timesteps causes the solution to quickly lose accuracy.
 
The procedure for the basis evolution approach is summarized in Algorithm~\ref{empirical_basic_evolution_algorithm}, which couples empirical chaos expansion with basis evolution in a finite-differencing method.

\begin{algorithm}
    \caption{\label{empirical_basic_evolution_algorithm}Empirical Chaos with Basis Evolution}
    \begin{algorithmic}[1]
\Function{ResampleBasis}{timestep}
\Comment{Return true if we should resample solutions to recompute the basis functions and return false if we should use the basis evolution operator. This can be a function that alternates between true and false, as in Figure~\ref{basis_evo_alt_tf50}.}
\If {$i$ is odd}
\State \Return True
\Else {}
\State \Return False
\EndIf
\EndFunction

\Procedure{EmpBasis}{}
\State $\xi[] \gets \text{vector with }n\ \text{values from domain of random variable}$
\State $\Delta t \gets \text{timestep size}$
\State $n_t \gets \text{number of timesteps}$
\For{$i = 1 \ldots n_t$}
\If{ResampleBasis(i)}
\For {$k = 1 \ldots n$}
\State $s[k] \gets \text{Solution of \eqref{ow_wave} for fixed } \xi[k]$
\EndFor
\State $T \gets \text{Matrix in \eqref{pod_sys} constructed from } s$
\State $U,\Sigma,V^t \gets \text{svd}(T)$
\State $k \gets \text{index of first scaled singular value smaller than } 10^{-4}$
\State $\Psi \gets \text{first } k-1 \text{ columns of } V^t$
\Else {}
\State $\text{Apply operator in \eqref{owwave_basis_exp_op_eq} to } \Psi$
\EndIf
\If {i = 1}
\State $\hat{u} \gets \text{ Projection of initial conditions onto } \Psi$
\Else {}
\State $\hat{u} \gets \text{ Projection of }\hat{u} \text{ onto } \Psi$
\EndIf
\State $\hat{u} \gets \text{ Solution of \eqref{ow_wave_sys}}$
\EndFor
\EndProcedure
\end{algorithmic}
\end{algorithm}

\section{Future Work}
Since their inception, polynomial chaos techniques have been successfully applied to a number of numerical problems that arise in uncertainty quantification. This work introduces a method to generate a set of empirical basis functions that varies with time to take the place of the standard orthogonal polynomial basis that remains fixed for the entire integration time. We demonstrate its numerical accuracy and efficiency over long-term integrations with a few model problems. We also introduce a method to numerically evolve the empirical basis functions without needing to resample solutions of the original SPDE. We demonstrate that the running time of the empirical chaos method scales linearly with the final integration time. Thus, it has the potential to address one of the two principal issues with applying polynomial chaos techniques--solving problems out to long-term integrations is often inefficient due to the need to continually increase the number of polynomial basis functions. 

This work does not, however, present a solution to the second and most fundamental issue--the curse of dimensionality. Namely, if the number of random variables is large then polynomial chaos techniques run slower than Monte Carlo methods. While empirical basis functions can be used in order to try to limit the size of such a basis, performing numerical integrations (which are a necessary component of the empirical chaos expansion) over a very high-dimensional space still poses significant challenges, both in terms of computational time and accuracy. We currently use sparse grid quadrature, which allows us to achieve accurate solutions over higher-dimensional spaces, but at the cost of sampling a large number of solutions of the original SPDE. Past a certain point, it becomes more practical to use Monte Carlo methods instead. Some techniques such as multi-element polynomial chaos have shown some promise when faced with larger numbers of random variables, and it is possible that using the basis decomposition employed by multi-element polynomial chaos methods coupled with empirical basis functions might allow higher-dimensional SPDEs to be solved more efficiently.

Another open question is to what degree the empirical chaos expansion algorithms can be made adaptive. Recall that the number of basis functions is selected by examining the scaled magnitude of the singular values from POD, and there are two values that we can adjust: (i) the length of each timestep; (ii) the number of sample trajectories of the SPDE. The timestep can be altered depending on the number of basis functions that we desire. In general, choosing a shorter timestep will result in the singular values from POD decaying faster, and thus will result in fewer empirical basis functions. The number of solutions to sample is a bit harder to determine. An approach that could be employed is to generate a set of empirical basis functions by sampling a fixed number of solutions to the SPDE, and then comparing the span of that basis to one that is generated by sampling a smaller number of solutions to the SPDE. If the difference in span is small then we would conclude that the number of solutions that we sampled is sufficient. If the difference in span is large, then we would conclude that the number of solutions that we sampled is insufficient, and proceed to sample additional solutions before repeating the process. It would be desirable to develop robust error estimates that could be used instead of generating two separate empirical bases.

The method of empirical basis expansion has the potential to be coupled with gPC expansions as well. We could use a set of $N$ orthogonal polynomial functions and perform the same trajectory sampling that we do for an empirical chaos expansion. We could then project the trajectories onto the polynomial basis, and construct an empirical basis for the residuals of the trajectories. The empirical basis should be numerically orthogonal to the polynomial basis, which might allow the two systems to be propagated independently. Even if the two sets of basis functions need to be coupled, this work has demonstrated that such systems can be solved by including a mass matrix that turns the system of propagation PDEs into a DAE.

Recent work \cite{pasini2013polynomial} has also shown that when gPC is applied to SPDEs with Hamiltonian structure, the resulting stochastic Galerkin system retains the Hamiltonian structure. Future work could attempt to determine whether the Hamiltonian structure is also retained when the original system is expanded using the empirical chaos method presented here. 

We have shown the empirical chaos method with basis evolution to be robust, accurate, and efficient for the model problems presented. Further work is needed to explore its efficacy for a broader range of problems. 

\bibliographystyle{siamplain}
\bibliography{references}

\end{document}